\documentclass[11pt,reqno,final]{article}

\usepackage{amsmath,amsfonts,amssymb,bm,amsthm,version}
\usepackage{mathrsfs,fancybox,pifont}
\usepackage{graphicx}
\usepackage{algorithm}
\usepackage{tikz}
\usetikzlibrary{calc, arrows.meta}
\usepackage{algpseudocode}
\usepackage{cite}
\usepackage{setspace}
\usepackage{url,hyperref}
\usepackage{fullpage} % set 1 inch margins
\usepackage[notcite,notref]{showkeys}
\usepackage{xcolor}
\usepackage{subfigure,multirow}
\usepackage{epstopdf}
\usepackage{subcaption}
\usepackage{cases}
\usepackage{geometry}
\usepackage{mathtools}
\usepackage{authblk}
\usepackage{lipsum}
\usepackage{indentfirst}
\newcommand{\dt}{\operatorname{d}\! t}

\newcommand{\dw}{\operatorname{d}\! W}

\allowdisplaybreaks

\numberwithin{equation}{section}
\theoremstyle{plain}
\newtheorem{theorem}{Theorem}%[section]
\newtheorem{proposition}{Proposition}%[section]
\newtheorem{assumption}{Assumption}%[section]

\theoremstyle{definition}
\newtheorem{definition}{Definition}[section]

\theoremstyle{remark}
\newtheorem{remark}{Remark}%[section]

\title{Constrained Mean--Variance Portfolio Selection under Reflected Factor Dynamics: Stochastic Riccati Equations, Neumann Problems, and a Deep Solver}

\author{
Zhecheng Huang, School of Mathematical Sciences Fudan University \& China Pacific Insurance (Group) Co., \href{mailto:huangzhecheng@fudan.edu.cn}{huangzhecheng@fudan.edu.cn}
\and
Guojiang Shao, School of Mathematical Sciences Fudan University, \href{mailto:gjshao23@m.fudan.edu.cn}{gjshao23@m.fudan.edu.cn}
\and   
Lei Wang, China Pacific Insurance (Group) Co.
\and Qi Zhang, School of Mathematical Sciences Fudan University \& State Key Laboratory of Cryptography and Digital Economy Security, Shandong University
}

\author{
Zhecheng Huang, School of Mathematical Sciences Fudan University \& China Pacific Insurance (Group) Co., \href{mailto:huangzhecheng@fudan.edu.cn}{huangzhecheng@fudan.edu.cn}
\and
Guojiang Shao, School of Mathematical Sciences Fudan University \& Department of Applied Mathematics The Hong Kong Polytechnic University, \href{mailto:gjshao23@m.fudan.edu.cn}{gjshao23@m.fudan.edu.cn} \thanks{Corresponding author}
\and Qi Zhang, School of Mathematical Sciences Fudan University \& State Key Laboratory of Cryptography and Digital Economy Security, Shandong University, \href{mailto:qzh@fudan.edu.cn}{qzh@fudan.edu.cn}
}

\begin{document}
\maketitle 

% \begin{abstract}
%     We study stochastic Riccati equations arising in constrained continuous-time mean--variance portfolio selection with random market coefficients. To model factor dynamics evolving within a bounded state space, we introduce reflected Markovian diffusions which gives rise naturally to homogeneous Neumann boundary conditions. A logarithmic transformation reformulates the extended stochastic Riccati equations as quadratic backward stochastic differential equations driven by the reflected factor processes. We establish a probabilistic representation linking these quadratic BSDEs to semilinear parabolic equations with Neumann boundary conditions and prove existence and uniqueness of bounded continuous viscosity solutions. Existence is obtained through stability and comparison arguments for quadratic BSDEs, whereas uniqueness is established by combining a Kobylanski-type transformation with a viscosity comparison principle adapted to the Neumann boundary setting. Building on this representation, we develop a Deep BSDE framework for the transformed equations. The method combines numerical simulation of reflected factor processes, a geometry-aware treatment of the Neumann boundary, and FISTA-based computation of the constrained Hamiltonian. Numerical experiments illustrate the accuracy and robustness of the propose approach. The resulting framework provides a mesh-free probabilistic method for high-dimensional Neumann problems arising in constrained portfolio selection and, more broadly, stochastic control. 
% \end{abstract}

\begin{abstract}

We study extended stochastic Riccati equations (ESREs) arising in constrained continuous-time mean--variance (MV) portfolio selection with random market coefficients. Reflected Markovian diffusions are introduced to model factor dynamics in a bounded state space, naturally leading to homogeneous Neumann boundary conditions. A logarithmic transformation converts the ESREs into quadratic backward stochastic differential equations (BSDEs) driven by the reflected factors. We establish the correspondence between these quadratic BSDEs and semilinear parabolic PDEs with homogeneous Neumann boundary conditions and prove existence and uniqueness of bounded continuous viscosity solutions. Existence is obtained through stability and comparison arguments for quadratic BSDEs, whereas uniqueness is established by combining a Kobylanski-type transformation with a viscosity comparison principle adapted to the Neumann boundary setting. We further develop a deep BSDE method combining reflected-path simulation, geometry-aware boundary treatment, and efficient evaluation of the constrained Hamiltonian. Numerical experiments illustrate the performance of the method, providing a mesh-free framework for high-dimensional stochastic Riccati equations arising in constrained portfolio optimization.

\end{abstract}


\section{Introduction}

Continuous-time mean--variance (MV) portfolio selection is a classical stochastic control problem balancing expected return and investment risk. Zhou and Li \cite{zhou2000continuous} showed that the pre-committed MV problem can be reformulated as a family of stochastic linear--quadratic (LQ) control problems, making Riccati equations a natural tool for characterizing optimal portfolios and the efficient frontier. When the market coefficients are random, the corresponding Riccati equations take the form of backward stochastic differential equations (BSDEs); see, e.g., Bismut \cite{bismut1976linear}. Portfolio constraints, such as short-selling restrictions and more general cone constraints, further complicate the problem and lead to extended stochastic Riccati equations (ESREs). Li, Zhou, and Lim \cite{li2002dynamic} studied mean--variance portfolio selection under no-shorting constraints, while Hu and Zhou \cite{hu2005constrained} developed a stochastic LQ framework with random coefficients and cone-constrained controls, in which the optimal strategy is characterized by a pair of nonlinear ESREs.

The ESREs arising in constrained portfolio optimization are analytically challenging due to their nonlinear generators, which may involve reciprocal dependence on the Riccati component and quadratic growth in the martingale variable. In addition, the Riccati solution must remain strictly positive. These features generally place such equations outside the standard Lipschitz BSDE framework and connect their analysis closely with the theory of BSDEs with quadratic growth. Important well-posedness results were established by Kohlmann and Tang \cite{kohlmann2002global}, Tang \cite{tang2003general}, Lim and Zhou \cite{lim2002mean}, Lim \cite{lim2004quadratic}, and Hu and Zhou \cite{hu2005constrained}, with quadratic-BSDE techniques such as comparison, truncation, and stability playing a central role; see also Kobylanski \cite{kobylanski2000backward}.

Most classical solvability results rely on boundedness or sufficiently strong integrability assumptions on the random coefficients. This becomes restrictive in factor-based market models, where interest rates, risk premia, and covariance matrices depend on stochastic factors that often evolve on unbounded state spaces. The resulting coefficients may therefore fail to remain uniformly bounded, making both the analytical treatment and numerical approximation of the associated SREs (ESREs) considerably more delicate. These difficulties motivate the reflected-factor formulation developed in this paper.

To address these difficulties, we adopt a reflected Markovian factor framework, in which the latent factors are constrained to evolve within a bounded admissible domain. The state constraint is enforced through a reflected stochastic differential equation (RSDE), which preserves the Markovian structure of the market model while keeping the factor process inside the prescribed state space. The well-posedness theory for reflected SDEs on constrained domains was established in the classical works of Lions and Sznitman \cite{lions1984stochastic} and Saisho \cite{saisho1987stochastic}. When the market coefficients are jointly continuous in time and the factor state, confining the factor process to a compact state space ensures that the resulting
market coefficients remain uniformly bounded. At the PDE level, the reflection mechanism naturally induces a homogeneous Neumann boundary condition. The probabilistic connection between reflected diffusions, BSDEs, and nonlinear Neumann boundary value problems is classical; see, for example, Pardoux and Zhang \cite{pardoux1998generalized} and Pardoux and R{\u{a}}{\c{s}}canu \cite{pardoux2014stochastic}. This reflected-factor formulation therefore provides a convenient Markovian setting, in which the associated SREs can be analyzed through quadratic BSDEs and semilinear Neumann problems.

Reflected diffusions provide a natural framework for modeling financial state variables constrained to economically meaningful regions. Unlike artificial truncation, reflection preserves the original drift--diffusion dynamics in the interior while regulating the process at the boundary. Such mechanisms have a long history in financial mathematics, dating back to El Karoui \cite{elkaroui2006processus}. Early applications include target-zone exchange-rate models \cite{krugman1991target,bertola1993stochastic} and regulated financial markets involving exchange rates, interest rates, and other prices \cite{bo2011some,bo2011conditional,bo2013conditional}. Related state constraints also arise in uncertain-volatility models with stochastic bounds \cite{fouque2018uncertain}.

% Reflection is particularly relevant to stochastic volatility modeling. For instance, the classical Stein--Stein model employs an Ornstein--Uhlenbeck (OU) volatility factor with Gaussian dynamics. When the factor is required to remain within a prescribed admissible range, a reflected OU process provides a natural alternative that preserves mean reversion in the interior of the domain. A reflected OU process preserves mean reversion while enforcing admissible bounds. 

Reflection is particularly relevant to stochastic volatility modeling. For instance, a reflected Ornstein--Uhlenbeck (OU) process provides a bounded-state alternative to the Gaussian volatility factor in the classical Stein--Stein model, preserving mean reversion while enforcing admissible bounds. Reflected stochastic-volatility models have been studied by Gulisashvili \cite{gulisashvili2023large}, while reflected and skew-diffusion structures have been applied to local-stochastic-volatility modeling and option pricing \cite{ding2021markov}. Extensions to multidimensional, path-dependent, and multisided reflected dynamics were developed by Ning and Wu \cite{ning2021well}. These studies motivate our use of reflected Markovian factors in portfolio optimization, where market coefficients depend on latent states constrained to economically or statistically admissible regions.

Deep learning methods for BSDEs provide a flexible probabilistic approach to high-dimensional nonlinear parabolic PDEs. The seminal deep BSDE method of E, Han, and Jentzen \cite{han2017deep,han2018solving} exploits the nonlinear Feynman--Kac correspondence between semilinear PDEs and forward--backward stochastic differential equations (FBSDEs). For Markovian problems, the BSDE solution provides a probabilistic representation of the PDE solution, while the martingale component contains information about its spatial gradient. Deep BSDE methods exploit this correspondence by approximating the martingale component with neural networks along simulated forward trajectories. Since the method samples the state space through stochastic trajectories rather than a spatial grid, it is particularly attractive for high-dimensional problems where conventional grid-based PDE solvers become computationally costly.

A number of related deep-learning formulations have subsequently been developed. Beck, E, and Jentzen \cite{beck2019machine} extended the probabilistic learning approach to fully nonlinear PDEs through second-order BSDE representations. Hur\'e, Pham, and Warin \cite{hure2020deep} proposed deep backward dynamic programming (DBDP) schemes, including the DBDP1 and DBDP2 algorithms, which approximate the solution recursively backward in time and provide convergence and error estimates. For stochastic optimal control problems, Ji, Peng, Peng, and Zhang \cite{ji2022solving} developed a deep learning method based on the stochastic maximum principle, reformulating the control problem through an extended Hamiltonian FBSDE system with a maximum condition.

Deep-learning methods based on probabilistic representations have also been extended to PDEs with Neumann boundary conditions. Boussange et al.~\cite{boussange2023deep} developed mesh-free methods for non-local nonlinear Neumann problems by combining reflected stochastic processes with deep splitting and multilevel Picard schemes. Deng and He~\cite{deng2026novel} extended FBSDE-based neural networks to Navier--Stokes and Cahn--Hilliard-type equations, employing reflection-based discretization for Neumann boundaries. More recently, Wang and Deng~\cite{wang2026deep} proposed a constraint-informed BSDE method for high-dimensional PDEs on manifolds and bounded domains with Neumann and Dirichlet boundary conditions. These developments demonstrate the potential of stochastic representations for high-dimensional PDEs with boundary constraints. In contrast, the present work focuses on constrained stochastic Riccati equations driven by reflected factors, where the logarithmic transformation leads to quadratic BSDEs. Our framework further combines geometry-aware boundary treatment with efficient evaluation of the constrained Hamiltonian, while accounting for positivity and quadratic growth.

The main contributions of this paper are threefold. First, we introduce reflected Markovian factors into constrained mean--variance portfolio selection with random coefficients. The reflection confines the factors to a bounded domain and induces homogeneous Neumann boundary conditions, while the logarithmic transformation reduces the ESREs to quadratic BSDEs suitable for analysis and computation. Second, we establish the BSDE--PDE correspondence for the transformed Riccati equations under reflected factor dynamics. We show that the Markovian decoupling field of the quadratic BSDE is a bounded continuous viscosity solution of the associated semilinear parabolic PDE with homogeneous Neumann boundary conditions, and prove uniqueness through a Kobylanski-type transformation and a viscosity comparison principle. This extends the reflected-diffusion BSDE--PDE framework to the quadratic-growth structure generated by stochastic Riccati equations. Finally, we develop a deep BSDE framework that combines reflected-path simulation, boundary-aware treatment of the Neumann condition, and efficient evaluation of the constrained Hamiltonian. The resulting method provides a mesh-free solver for high-dimensional stochastic Riccati equations.

The rest of the paper is organized as follows. Section~2 formulates the constrained mean--variance portfolio
problem in a market with random coefficients and presents the associated pair of ESREs.
Section~3 introduces reflected Markovian factors on a bounded domain, transforms the ESREs into quadratic BSDEs, and establishes their connection to semilinear parabolic
PDEs with homogeneous Neumann boundary conditions, including the existence and uniqueness of bounded continuous viscosity solutions. Section~4 develops a deep BSDE method for the transformed equations and presents numerical experiments. Section~5 concludes the paper and outlines directions
for future research. The proofs of the transformation theorem and the PDE
comparison principle are deferred to Appendices~\ref{app-kobi} and~\ref{app-comparison}, respectively.

\section{Problem Formulation}

Let $(\Omega,\mathcal F,\mathbb P)$ be a complete
probability space supporting an $n$-dimensional
standard Brownian motion $W=(W(t))_{t\in[0,T]}$.
Let $\mathbb F=(\mathcal F_t)_{t\in[0,T]}$
denote the augmented natural filtration of $W$. The superscript $^*$ denotes the transpose of a
vector or matrix. For a symmetric matrix $A$,
the notation $A\geq 0$ (resp., $A>0$) means
that $A$ is positive semidefinite
(resp., positive definite).
We denote by $I_n$ the $n\times n$ identity matrix,
and by $\langle\cdot,\cdot\rangle$ and $|\cdot|$
the Euclidean inner product and norm, respectively.
Let $\mathbb S^n$ be the space of real symmetric
$n\times n$ matrices, and let $C_b^2(\mathbb R^n)$
denote the space of twice continuously differentiable
functions whose derivatives up to order two,
including the functions themselves, are bounded. For a continuous local martingale $\mathfrak m$ with $\mathfrak m_0=0$, we denote by $\mathcal E(\mathfrak m)$ its stochastic exponential, i.e.,
\[
    \mathcal E(\mathfrak m)_t
    :=\exp\left(\mathfrak m_t-\frac12\langle \mathfrak m\rangle_t\right),
    \qquad t\in[0,T],
\]
where $\langle \mathfrak m\rangle$ denotes the quadratic variation of $\mathfrak m$. For any $k,j\in\mathbb N$, we use the following spaces throughout the paper:
\begin{align*}
\begin{array}{rl} 
% L_{\mathcal{F}_T}^{\infty}\left( \mathbb{R}^m \right): &\hspace{-0.25cm} \text {the set of $\mathcal{F}_T$-measurable and essentially bounded } \mathbb{R}^m \text{-valued random variables;}\\ 
% L_{\mathbb{F}}^{\infty}\left(0, T ; \mathbb{R}^m \right): &\hspace{-0.25cm} \text {the set of $\mathbb{F}$-adapted essentially bounded } \mathbb{R}^m \text{-valued processes;}\\ 
L_{\mathbb{F}}^{\infty}\left(0, T ; \mathbb{R}^{k \times j} \right): &\hspace{-0.25cm} \text {the set of $\mathbb{F}$-progressively measurable, bounded and } \mathbb{R}^{k \times j} \text{-valued processes};\\
% L_{\mathbb{F}}^{\infty}\left(0, T ; \mathbb{R}_{+}\right): &\hspace{-0.25cm} \text {the set of $\mathbb{F}$-adapted essentially bounded nonnegative processes};\\
% L_{\mathbb{F}}^{\infty}\left(0, T ; \mathbb{R}_{>0}\right): &\hspace{-0.25cm} \text {the set of $\mathbb{F}$-adapted essentially bounded positive processes};\\
% L_{\mathbb{F}}^{\infty}\left(0, T ; \mathbb{R}_{\geq 0}\right): &\hspace{-0.25cm} \text {the set of $\mathbb{F}$-adapted essentially bounded non-negative processes};\\
L_{\mathbb{F}}^{\infty}\left(0, T ; \mathbb{R}_{\gg 1}\right): &\hspace{-0.25cm} \text {the set of $\mathbb{F}$-progressively measurable processes } v:[0, T] \times \Omega \rightarrow(0,+\infty) \\ \hspace{-0.3cm} &\hspace{-0.25cm}  ~\text{such that }
 \text {$c^{-1} \leqslant v(t) \leqslant c$ a.e. a.s. for some constant $c>0$};\\
L_{\mathbb{F}}^2\left(0, T ; \mathbb{R}^k \right): &\hspace{-0.25cm} \text {the set of $\mathbb{F}$-adapted processes } v:[0, T] \times \Omega \rightarrow \mathbb{R}^k ~\text {such that } \\
&\hspace{-0.25cm}\mathbb{E}\left[\int_0^T|v(t)|^2 \dt\right]<\infty;\\
% L_{\mathbb{F}}^1\left(0, T ; \mathbb{R}^m \right): &\hspace{-0.25cm} \text {the set of $\mathbb{F}$-adapted processes } v:[0, T] \times \Omega \rightarrow \mathbb{R}^m ~\text {such that } \\
% &\hspace{-0.25cm}\mathbb{E}\left[\int_0^T|v(t)|\dt\right]<\infty;\\
% L_{\mathbb{F}}^{2,c}\left(0,T; \mathbb{R}\right): &\hspace{-0.25cm} \text {the set of $\mathbb{F}$-adapted processes $v:[0, T] \times \Omega \rightarrow \mathbb{R}$ with} \\
% &\hspace{-0.25cm}\text {continuous sample paths such that } \mathbb{E}\left[\sup_{t \in[0, T]}|v(t)|^2\right]<\infty; \\
% L_{\mathbb{F}}^{\infty,c}\left(0,T; \mathbb{R}\right): &\hspace{-0.25cm} \text {the set of $\mathbb{F}$-adapted essentially bounded  processes $v:[0, T] \times \Omega \rightarrow \mathbb{R}$ with} \\
% &\hspace{-0.25cm}\text {continuous sample paths}. 
\end{array}
\end{align*}

We consider a financial market consisting of one risk-free bond and $m$ risky stocks over a finite horizon $[0,T]$. The price of the risk-free asset, denoted by $S_0(t)$, evolves according to
\[
\left\{
\begin{aligned}
\mathrm{d}S_0(t) &= r(t)S_0(t)\,\mathrm{d}t, \qquad t\in(0,T],\\
S_0(0) &= s_0>0,
\end{aligned}
\right.
\]
where $r(\cdot)\in L_{\mathbb F}^{\infty}(0,T;\mathbb R)$ denotes the short rate. For $i=1,\ldots,m$, the price process of the $i$-th risky asset satisfies
\[
\left\{
\begin{aligned}
\mathrm{d}S_i(t)
&=
S_i(t)
\left\{
\mu_i(t)\,\mathrm{d}t
+
\sum_{j=1}^{n}\sigma_{ij}(t)\,\mathrm{d}W_j(t)
\right\},
\qquad t\in(0,T],\\
S_i(0)&=s_i>0,
\end{aligned}
\right.
\]
where $\mu_i(t)$ denotes the appreciation rate and
\[
\sigma_i(t) := 
\bigl(\sigma_{i1}(t),\ldots,\sigma_{in}(t)\bigr)
\]
is the volatility vector of the $i$-th risky asset. We assume that the market coefficients $r$, $\mu_i$, and $\sigma_{ij}$ are bounded and $\mathbb F$-adapted.

\begin{remark}
The boundedness assumption on the market coefficients is used in the well-posedness analysis below. No sign restriction on $r(t)$, $\mu_i(t)$, or $\sigma_{ij}(t)$ is essential to our main results.
\end{remark}

Define the appreciation-rate vector and volatility matrix by
\[
\mu(t)
:=
\begin{bmatrix}
\mu_1(t)\\
\vdots\\
\mu_m(t)
\end{bmatrix}
\in\mathbb{R}^m,
\qquad
\sigma(t)
:=
\begin{bmatrix}
\sigma_1(t)\\
\vdots\\
\sigma_m(t)
\end{bmatrix}
\in\mathbb{R}^{m\times n}.
\]
Then
\[
\mu(\cdot)\in L_{\mathbb F}^{\infty}(0,T;\mathbb R^m),
\qquad
\sigma(\cdot)\in L_{\mathbb F}^{\infty}(0,T;\mathbb R^{m\times n}).
\]
We assume that $\sigma$ is uniformly nondegenerate: there exists a constant $\varepsilon>0$ such that
\[
\sigma(t)\sigma^*(t)\geq \varepsilon I_m,
\qquad
t\in[0,T],\quad \mathbb P\text{-a.s.}
\]

Let $X(t)$ denote the investor's total wealth at time $t$. Under continuous trading, with no transaction costs or intermediate consumption, the self-financing wealth process satisfies
\begin{equation}\label{state}
\left\{
\begin{aligned}
\mathrm{d}X(t)
&=
\left[
r(t)X(t)
+
\sum_{i=1}^m \bigl(\mu_i(t)-r(t)\bigr)\pi_i(t)
\right]\mathrm{d}t \\
&\quad
+
\sum_{i=1}^m\sum_{j=1}^n
\sigma_{ij}(t)\pi_i(t)\,\mathrm{d}W_j(t),
\qquad t\in(0,T],\\
X(0)&=x_0>0,
\end{aligned}
\right.
\end{equation}
where $\pi_i(t)$ denotes the amount invested in the $i$-th risky asset. A negative value of $\pi_i(t)$ corresponds to a short position. The portfolio strategy is denoted by
\[
\pi(t):=
\bigl(\pi_1(t),\ldots,\pi_m(t)\bigr)^*.
\]
The amount invested in the risk-free asset is determined by $\pi_0(t)=X(t)-\sum_{i=1}^m\pi_i(t)$. Let $\mathbf{1}_m=(1,\ldots,1)^*\in\mathbb R^m$ and define the excess-return vector
\[
\widetilde{\mu}(t):=\mu(t)-r(t)\mathbf{1}_m.
\]
Then \eqref{state} can be written compactly as
\begin{equation}\label{dynamics}
\left\{
\begin{aligned}
\mathrm{d}X(t)
&=
\bigl[r(t)X(t)+\langle\widetilde{\mu}(t),\pi(t)\rangle\bigr]\mathrm{d}t
+
\langle\sigma^*(t)\pi(t),\mathrm{d}W(t)\rangle,\\
X(0)&=x_0.
\end{aligned}
\right.
\end{equation}

Let $\Gamma\subset\mathbb R^m$ be a nonempty closed cone,
that is, $\Gamma$ is closed and $\lambda\pi\in\Gamma$
for every $\pi\in\Gamma$ and $\lambda\geq 0$.
A typical example is $\Gamma=\mathbb R_+^m$. The set of admissible portfolio strategies is defined by
\[
\mathcal{A}
:=
\left\{
\pi(\cdot)\in L_{\mathbb F}^2(0,T;\mathbb R^m)
\;\middle|\;
\pi(t)\in\Gamma \ \text{a.e. a.s., and \eqref{dynamics} admits a unique strong solution}
\right\}.
\]
For $\pi(\cdot)\in\mathcal A$, let $X(\cdot)$ denote the corresponding solution of \eqref{dynamics}; then $(X(\cdot),\pi(\cdot))$ is called an admissible pair.

Given a target expected terminal wealth $d\in\mathbb R$, the investor seeks an admissible portfolio that minimizes the variance of terminal wealth subject to $\mathbb E[X(T)]=d$. The constrained mean--variance problem is therefore
\begin{equation}\label{prob1}
\begin{aligned}
\min_{\pi(\cdot)\in\mathcal A}\quad
&\mathcal J(\pi(\cdot),d)
:=
\operatorname{Var}(X(T))
=
\mathbb E\bigl[|X(T)-d|^2\bigr],\\
\text{subject to}\quad
&\mathbb E[X(T)]=d.
\end{aligned}
\end{equation}

{ \bf Problem \eqref{prob1} } is a constrained stochastic control problem, whose feasibility and solution have been studied extensively in the literature (see, e.g., Hu and Zhou \cite{hu2005constrained}, Hu, Shi and Xu \cite{hu2022constrained}). In particular, the feedback structure of the optimal strategy and the value function of the problem are associated with the following two extended stochastic Riccati equations with random coefficients:
\begin{equation} \label{SRE-positive}
    \left\{\begin{aligned}
        \mathrm{d}P_+ &= - \left[2rP_+ + \inf_{\pi \in \Gamma}\Big(\langle P_+ \sigma \sigma^* \pi,\pi \rangle + 2 \langle P_+ \tilde{\mu} + \sigma \Lambda_+,\pi \rangle\Big) \right]\dt + \langle  \Lambda_+,\dw\rangle, \quad t \in [0, T),  \\
        P_+(T) & = 1, \qquad P_+(t) > 0, \quad t \in [0, T), \\
    \end{aligned} \right.
\end{equation}
and
\begin{equation} \label{SRE-negative}
    \left\{\begin{aligned}
         \mathrm{d}P_- &= - \left[2rP_- + \inf_{\pi \in \Gamma}\Big(\langle P_- \sigma \sigma^* \pi,\pi \rangle - 2 \langle P_- \tilde{\mu} + \sigma \Lambda_-,\pi \rangle\Big) \right]\dt + \langle \Lambda_-,\dw\rangle, \quad t \in [0, T),  \\
        P_-(T) & = 1, \qquad P_-(t) > 0, \quad t \in [0, T). \\
    \end{aligned} \right.
\end{equation}

For notational simplicity, we suppress the time arguments of processes and coefficients whenever
no confusion arises.
The ESREs \eqref{SRE-positive} and \eqref{SRE-negative} are difficult to analyze and solve numerically due to the nonlinear dependence
of their generators on the Riccati processes and martingale components. These difficulties motivate the logarithmic
transformation and reflected-factor framework
developed in the next section.

\section{Reflected Markovian Factors, Quadratic BSDEs, and Neumann PDEs}

\subsection{Reflected Markovian Factor Process}

A common approach to random-coefficient market models is to impose a Markovian factor structure, under which the randomness of the market coefficients is generated by an underlying forward diffusion. This representation also provides a natural link between the associated stochastic control problem and PDE methods; see, for example, Section~6.4 of Kohlmann and Tang \cite{kohlmann2002global} and Section~7 of Lim \cite{lim2004quadratic}.

In the standard Markovian setting, for each $(t,v)\in[0,T]\times\mathbb R^n$, the factor process $(\bm V_s^{t,v})_{s\in[t,T]}$ satisfies
\[
\mathrm{d}\bm V_s^{t,v}
=
\bm F(\bm V_s^{t,v})\,\mathrm{d}s
+
\bm G(\bm V_s^{t,v})\,\mathrm{d}W_s,
\qquad
\bm V_t^{t,v}=v,
\]
where
\[
\bm F:\mathbb R^n\to\mathbb R^n,
\qquad
\bm G:\mathbb R^n\to\mathbb R^{n\times n}
\]
are deterministic functions. The market coefficients are then represented as deterministic functions of time and the current factor state:
\[
r_s=r(s,\bm V_s^{t,v}),\qquad
\widetilde{\mu}_s=\widetilde{\mu}(s,\bm V_s^{t,v}),\qquad
\sigma_s=\sigma(s,\bm V_s^{t,v}),
\]
where
\[
r:[0,T]\times\mathbb R^n\to\mathbb R,\qquad
\widetilde{\mu}:[0,T]\times\mathbb R^n\to\mathbb R^m,\qquad
\sigma:[0,T]\times\mathbb R^n\to\mathbb R^{m\times n}.
\]
Let $\mathcal D\subset\mathbb R^n$ be a bounded open domain whose geometric properties will be specified below.
\begin{assumption} \label{ass-market}
Assume that $r$, $\widetilde{\mu}$, and $\sigma$ are bounded and Lipschitz continuous on $[0,T]\times\overline{\mathcal D}$. Moreover, suppose that $\sigma$ is uniformly nondegenerate, i.e.,
\[
\sigma(t,v)\sigma(t,v)^* \geq \delta I_m,
\qquad
(t,v)\in[0,T]\times\overline{\mathcal D},
\]
for some constant $\delta>0$.
\end{assumption}

The well-posedness of SREs (ESREs) with genuinely unbounded random coefficients remains technically challenging in general. Existing results typically rely on additional structural or integrability conditions; see, for example, Shen \cite{shen2015mean}, Shen and Zeng \cite{shen2015optimal}, Algoulity \cite{algoulity2021quadratic}, and Gnoatto et al.~\cite{gnoatto2025deep}.

Motivated by reflected diffusion models in finance, we instead let the factor process evolve on a bounded admissible domain through a reflected SDE. The domain represents economically relevant ranges of market states, while the reflection acts only at the boundary to keep the factor process inside the prescribed region. Consequently, the state-dependent market coefficients remain uniformly bounded on the relevant state space, placing the associated ESREs within a tractable quadratic-BSDE framework. At the same time, the reflection mechanism naturally induces a homogeneous Neumann boundary condition, providing the probabilistic link between the reflected forward process, the quadratic BSDE, and the corresponding PDE on $\overline{\mathcal D}$.

To formalize the reflection mechanism that keeps the process within a bounded domain, we impose the following geometric and analytical conditions on the state domain and the model coefficients:

\begin{assumption}[Domain Regularity]\label{ass:domain}
Let $\mathcal D\subset\mathbb R^n$ be a bounded open domain of the form
\[
\mathcal D=\{x\in\mathbb R^n:\phi(x)>0\},
\qquad
\partial\mathcal D=\{x\in\mathbb R^n:\phi(x)=0\},
\]
for some $\phi\in C_b^2(\mathbb R^n)$ satisfying
\[
|\nabla\phi(x)|=1,
\qquad x\in\partial\mathcal D.
\]
Thus, $\nabla\phi(x)$ is the inward unit normal vector at each $x\in\partial\mathcal D$.
\end{assumption}

\begin{assumption}[Coefficient Regularity]\label{ass:coefficients}
The drift $\bm F:\overline{\mathcal D}\to\mathbb R^n$ and diffusion coefficient
$\bm G:\overline{\mathcal D}\to\mathbb R^{n\times n}$ are Lipschitz continuous. That is, there exists a constant $C>0$ such that
\[
|\bm F(v)-\bm F(v')|
+
|\bm G(v)-\bm G(v')|
\leq C|v-v'|,
\qquad
v,v'\in\overline{\mathcal D}.
\]
\end{assumption}

Under Assumptions~\ref{ass:domain} and \ref{ass:coefficients}, for each
$(t,v)\in[0,T]\times\overline{\mathcal D}$, consider the reflected SDE
\begin{equation}\label{eq:RSDE}
\left\{
\begin{aligned}
\mathrm{d}\bm V_s^{t,v}
&=
\bm F(\bm V_s^{t,v})\,\mathrm{d}s
+
\bm G(\bm V_s^{t,v})\,\mathrm{d}W_s
+
\nabla\phi(\bm V_s^{t,v})\,\mathrm{d}\zeta_s^{t,v},
\qquad s\in[t,T],\\
\bm V_t^{t,v}&=v.
\end{aligned}
\right.
\end{equation}
By the classical theory of reflected SDEs developed by Lions and Sznitman \cite{lions1984stochastic} and Saisho \cite{saisho1987stochastic}, \eqref{eq:RSDE} admits a unique strong solution $(\bm V^{t,v},\zeta^{t,v})$. Moreover, $V^{t,v}$ has continuous sample paths satisfying $V_s^{t,v}\in\overline{\mathcal D}$ for all $s\in[t,T]$, almost surely, and $\zeta^{t,v}$ is a continuous nondecreasing adapted process with $\zeta_t^{t,v}=0$ whose increments are supported on the boundary:
\begin{equation}\label{eq:flatness} 
    \zeta_s^{t,v} = \int_t^s \mathbf 1_{\{\bm V_\varsigma^{t,v}\in \partial\mathcal D\}} \,\mathrm{d}\zeta_\varsigma^{t,v}, \qquad s\in[t,T]. 
\end{equation}
Hence, the reflection term acts only when the factor process reaches
$\partial\mathcal D$.

\subsection{Logarithmic Transformation of ESREs}

The well-posedness of ESREs \eqref{SRE-positive} and \eqref{SRE-negative} follows from existing results; see, for example, Theorem $3.6$ in Hu et al. \cite{hu2022constrained}.
\begin{theorem}
    Under Assumption~\ref{ass-market}, the ESREs \eqref{SRE-positive} and \eqref{SRE-negative} admit unique solutions $(P_+,\Lambda_+)$ and
$(P_-,\Lambda_-)$, respectively, in
\[
L^\infty_{\mathbb F}(0,T;\mathbb R_{\gg 1})
\times L^2_{\mathbb F}(0,T;\mathbb R^n).
\]
\end{theorem}

Since $P_+$ and $P_-$ are uniformly positive, we introduce the logarithmic transformation
\[
\check{P}_+(t):=\log P_+(t),\qquad
\check{\Lambda}_+(t):=\frac{\Lambda_+(t)}{P_+(t)},
\]
and similarly,
\[
\check{P}_-(t):=\log P_-(t),\qquad
\check{\Lambda}_-(t):=\frac{\Lambda_-(t)}{P_-(t)}.
\]
Applying It\^o's formula yields
\begin{equation}\label{log-positive}
\left\{
\begin{aligned}
\mathrm{d}\check{P}_+
&=
-\left[
2r
+
\inf_{\pi\in\Gamma}
\Big(
\langle\sigma\sigma^*\pi,\pi\rangle
+
2\langle\widetilde{\mu}+\sigma\check{\Lambda}_+,\pi\rangle
\Big)
+
\frac12|\check{\Lambda}_+|^2
\right]\mathrm{d}t
+
\langle\check{\Lambda}_+,\mathrm{d}W\rangle,
\\
\check{P}_+(T)&=0,
\end{aligned}
\right.
\end{equation}
and
\begin{equation}\label{log-negative}
\left\{
\begin{aligned}
\mathrm{d}\check{P}_-
&=
-\left[
2r
+
\inf_{\pi\in\Gamma}
\Big(
\langle\sigma\sigma^*\pi,\pi\rangle
-
2\langle\widetilde{\mu}+\sigma\check{\Lambda}_-,\pi\rangle
\Big)
+
\frac12|\check{\Lambda}_-|^2
\right]\mathrm{d}t
+
\langle\check{\Lambda}_-,\mathrm{d}W\rangle,
\\
\check{P}_-(T)&=0.
\end{aligned}
\right.
\end{equation}

Thus, the original ESREs \eqref{SRE-positive}--\eqref{SRE-negative} are reformulated as quadratic BSDEs. Conversely, any bounded solution $\check P_\pm$
with square-integrable martingale component
$\check\Lambda_\pm$ yields a solution of the
original ESRE through
\[
P_\pm=e^{\check P_\pm},
\qquad
\Lambda_\pm=P_\pm\check\Lambda_\pm.
\]

The logarithmic transformation incorporates the positivity constraint directly into the parametrization and reveals the quadratic structure of the martingale component without singular dependence on the Riccati variable. This formulation is therefore more convenient for both theoretical analysis and numerical approximation. In what follows, we focus on \eqref{log-positive} under the reflected factor dynamics \eqref{eq:RSDE}; the analysis of \eqref{log-negative} is analogous.

\subsection{Neumann Boundary PDE Associated with the Reflected Factor Process}

For $(t,v,z)\in[0,T]\times\overline{\mathcal D}\times\mathbb R^n$, define the generator
\[
\mathcal H(t,v,z)
:=
2r(t,v)
+
\inf_{\pi\in\Gamma}
\left\{
\left\langle
\sigma(t,v)\sigma(t,v)^*\pi,\pi
\right\rangle
+
2\left\langle
\widetilde{\mu}(t,v)+\sigma(t,v)z,\pi
\right\rangle
\right\}
+
\frac12|z|^2.
\]
The following proposition summarizes the growth and continuity properties of $\mathcal H$.

\begin{proposition}\label{prop:cone-generator}
There exists a constant $C>0$ such that, for all
$(t_i,v_i,z_i)\in[0,T]\times\overline{\mathcal D}\times\mathbb R^n$,
$i=1,2$,
\[
|\mathcal H(t_1,v_1,z_1)|
\leq
C\bigl(1+|z_1|^2\bigr),
\]
and
\[
\begin{aligned}
& \quad
|\mathcal H(t_1,v_1,z_1)-\mathcal H(t_2,v_2,z_2)|
\\
&\leq
C\bigl(1+|z_1|^2+|z_2|^2\bigr)
\bigl(|t_1-t_2|+|v_1-v_2|\bigr)
+
C\bigl(1+|z_1|+|z_2|\bigr)|z_1-z_2|.
\end{aligned}
\]
\end{proposition}

\begin{proof}
    For fixed \((t,v,z)\), the map 
    \begin{equation} \label{mapping}
        \pi\mapsto \left\langle \sigma(t,v)\sigma(t,v)^*\pi,\pi \right\rangle + 2\left\langle \widetilde\mu(t,v)+\sigma(t,v)z,\pi \right\rangle,
    \end{equation}
    is uniformly convex because \(\sigma\sigma^*\) is uniformly positive definite. Since \(0\in\Gamma\), its infimum over \(\Gamma\) is non-positive. On the other hand, by uniform ellipticity, 
    \begin{equation*}
        \left\langle \sigma\sigma^*\pi,\pi \right\rangle + 2\left\langle \widetilde\mu+\sigma z,\pi \right\rangle \geq \delta|\pi|^2 - C(1+|z|)|\pi| \geq -C(1+|z|^2).
    \end{equation*}
    Hence the infimum is bounded in absolute value by \(C(1+|z|^2)\), which gives 
    \[ 
        |\mathcal H(t_1,v_1,z_1)| \leq C(1+|z_1|^2). 
    \] 
    For the continuity estimate, let
\[
    A_i:=\sigma(t_i,v_i)\sigma(t_i,v_i)^*, \qquad
    b_i:=\widetilde\mu(t_i,v_i)+\sigma(t_i,v_i)z_i, \qquad i=1,2.
\]
Denote
\begin{equation} \label{minimizer}
    \ell_i := \inf_{\pi\in\Gamma} \left\{ \langle A_i\pi,\pi\rangle + 2\langle b_i,\pi\rangle
    \right\}.
\end{equation}
Since $\Gamma$ is closed and the objective function is
coercive, the infimum in \eqref{minimizer} is attained.
Let $\pi_i\in\Gamma$ be a minimizer of \eqref{minimizer}. As shown above, the coercivity of \(A_i\) and the fact that \(0\in\Gamma\) imply
\[
    |\pi_i| \leq C(1+|z_i|), \qquad i=1,2.
\]
We estimate \(|\ell_1-\ell_2|\). By the optimality of \(\pi_1\),
\[
    \ell_1 = \langle A_1\pi_1,\pi_1\rangle + 2\langle b_1,\pi_1\rangle \leq \langle A_1\pi_2,\pi_2\rangle
    + 2\langle b_1,\pi_2\rangle.
\]
Therefore,
\[
\begin{aligned}
    \ell_1-\ell_2 &\leq \langle A_1\pi_2,\pi_2\rangle + 2\langle b_1,\pi_2\rangle - \left( \langle A_2\pi_2,\pi_2\rangle
    + 2\langle b_2,\pi_2\rangle \right)  \\
    &= \langle (A_1-A_2)\pi_2,\pi_2\rangle + 2\langle b_1-b_2,\pi_2\rangle.
\end{aligned}
\]
Similarly, by testing the second infimum at \(\pi_1\), we get
\[
    \ell_2-\ell_1 \leq \langle (A_2-A_1)\pi_1,\pi_1\rangle + 2\langle b_2-b_1,\pi_1\rangle.
\]
Hence
\[
    |\ell_1-\ell_2| \leq |A_1-A_2| \left(|\pi_1|^2+|\pi_2|^2\right) +2 |b_1-b_2| \left(|\pi_1|+|\pi_2|\right).
\]
By the Lipschitz continuity and boundedness of \(\sigma\) on \([0,T]\times\overline{\mathcal D}\),
\[
\begin{aligned}
    |A_1-A_2| &= |\sigma(t_1,v_1)\sigma(t_1,v_1)^* - \sigma(t_2,v_2)\sigma(t_2,v_2)^*| \\
    &\leq C\bigl(|t_1-t_2|+|v_1-v_2|\bigr).
\end{aligned}
\]
Moreover,
\begin{equation*}
    b_1-b_2 = \widetilde\mu(t_1,v_1)-\widetilde\mu(t_2,v_2)+ \bigl(\sigma(t_1,v_1)-\sigma(t_2,v_2)\bigr)z_1 + \sigma(t_2,v_2)(z_1-z_2).
\end{equation*}
Thus
\[
    |b_1-b_2| \leq C\bigl(1+|z_1|\bigr) \bigl(|t_1-t_2|+|v_1-v_2|\bigr) + C|z_1-z_2|.
\]
Consequently, we obtain the symmetric bound
\[
    |b_1-b_2| \leq C\bigl(1+|z_1|+|z_2|\bigr) \bigl(|t_1-t_2|+|v_1-v_2|\bigr) + C|z_1-z_2|.
\]
Combining this with
\[
    |\pi_i|\leq C(1+|z_i|), \qquad i=1,2,
\]
we obtain
\[
\begin{aligned}
    |\ell_1-\ell_2| \leq C\bigl(1+|z_1|^2+|z_2|^2\bigr)
    \bigl(|t_1-t_2|+|v_1-v_2|\bigr) + C(1+|z_1|+|z_2|)|z_1-z_2|.
\end{aligned}
\]
Finally,
\[
    |r(t_1,v_1)-r(t_2,v_2)| \leq C\bigl(|t_1-t_2|+|v_1-v_2|\bigr),
\]
and
\[
    \left| \frac12|z_1|^2-\frac12|z_2|^2 \right| \leq \frac12(|z_1|+|z_2|)|z_1-z_2|.
\]
Putting these estimates together gives
\[
\begin{aligned}
&
|\mathcal H(t_1,v_1,z_1) - \mathcal H(t_2,v_2,z_2)|
\\ \qquad\leq
&C\bigl(1+|z_1|^2+|z_2|^2\bigr) \bigl(|t_1-t_2|+|v_1-v_2|\bigr)
+ C(1+|z_1|+|z_2|)|z_1-z_2|.
\end{aligned}
\]
This completes the proof.
\end{proof}

\begin{remark}
The estimates in Proposition~\ref{prop:cone-generator} will be used in the viscosity analysis below. In particular, they show that the generator $\mathcal H$ satisfies the quadratic-growth and continuity conditions analogous to condition (H5) in Kobylanski \cite{kobylanski2000backward}. However, the viscosity-solution uniqueness result in \cite{kobylanski2000backward} does not apply directly in our setting, since the reflected factor process leads to a PDE posed on a bounded domain with a homogeneous Neumann boundary condition. A separate comparison argument adapted to the Neumann problem is therefore required.
\end{remark}

For technical simplicity, we study the connection between the logarithmically transformed BSDE \eqref{log-positive} and its associated PDE, rather than working directly with the original ESRE \eqref{SRE-positive}. Let
\[
\mathcal L\varphi(v)
:=
\frac12
\operatorname{Tr}
\left(
\bm G(v)\bm G(v)^*D^2\varphi(v)
\right)
+
\langle
\bm F(v),\nabla\varphi(v)
\rangle
\]
denote the infinitesimal generator of the reflected diffusion \eqref{eq:RSDE}.

For each $(t,v)\in[0,T]\times\overline{\mathcal D}$, let
\((Y_s^{t,v},Z_s^{t,v})_{s\in[t,T]}\)
denote the unique bounded solution of the quadratic BSDE
\begin{equation}\label{eq:flow_quadratic}
Y_s^{t,v}
=
\int_s^T
\mathcal H(\varsigma,\bm V_\varsigma^{t,v},Z_\varsigma^{t,v})
\,\mathrm d\varsigma
-
\int_s^T
\langle Z_\varsigma^{t,v},\mathrm dW_\varsigma\rangle,
\qquad s\in[t,T],
\end{equation}
with terminal condition
\[
Y_T^{t,v}=0.
\]

The relationship between this quadratic BSDE and the associated semilinear PDE with homogeneous Neumann boundary condition is established through the Markovian representation of the solution.
\begin{theorem}\label{thm:pde-bsde}
Assume that there exists a function
$u\in C^{1,2}([0,T]\times\overline{\mathcal D})$
satisfying
\begin{equation}\label{eq:pde-sre}
\left\{
\begin{aligned}
\partial_t u(t,v)
+\mathcal L u(t,v)
+\mathcal H
\bigl(t,v,\bm G(v)^*\nabla u(t,v)\bigr)
&=0,
&& (t,v)\in[0,T)\times\mathcal D,
\\
\langle\nabla u(t,v),\nabla\phi(v)\rangle
&=0,
&& (t,v)\in[0,T)\times\partial\mathcal D,
\\
u(T,v)&=0,
&& v\in\overline{\mathcal D}.
\end{aligned}
\right.
\end{equation}

For each $(t,v)\in[0,T]\times\overline{\mathcal D}$, let
$(\bm V^{t,v},\zeta^{t,v})$ be the solution of the reflected SDE
\eqref{eq:RSDE}. Define
\[
Y_s^{t,v}:=
u(s,\bm V_s^{t,v}),
\qquad
Z_s^{t,v}:=
\bm G(\bm V_s^{t,v})^*
\nabla u(s,\bm V_s^{t,v}).
\]
Then $(Y^{t,v},Z^{t,v})$ is an adapted solution of the quadratic BSDE
\eqref{eq:flow_quadratic}.
\end{theorem}

\begin{proof}
Fix $(t,v)\in[0,T]\times\overline{\mathcal D}$. Since
$u\in C^{1,2}([0,T]\times\overline{\mathcal D})$ and
$\overline{\mathcal D}$ is compact, both $u$ and $\nabla u$ are bounded.
Together with the continuity of $\bm G$ on $\overline{\mathcal D}$, this implies
that
\[
Z_s^{t,v}
=
\bm G(\bm V_s^{t,v})^*
\nabla u(s,\bm V_s^{t,v})
\]
is bounded, and hence
\[
\mathbb E\int_t^T |Z_s^{t,v}|^2\,\mathrm ds<\infty .
\]

Applying It\^o's formula to
$u(s,\bm V_s^{t,v})$ over the interval $[s,T]$ yields
\begin{equation} \label{eq:Ito}
    \begin{aligned}
        Y_T^{t,v}-Y_s^{t,v}
={}&
\int_s^T
\left(
\partial_\varsigma u+\mathcal L u
\right)
(\varsigma,\bm V_\varsigma^{t,v})
\,\mathrm d\varsigma + 
\int_s^T
\left\langle
Z_\varsigma^{t,v},
\mathrm dW_\varsigma
\right\rangle
\\
&+
\int_s^T
\left\langle
\nabla u(\varsigma,\bm V_\varsigma^{t,v}),
\nabla\phi(\bm V_\varsigma^{t,v})
\right\rangle
\,\mathrm d\zeta_\varsigma^{t,v}.
    \end{aligned}
\end{equation}

By the homogeneous Neumann boundary condition in \eqref{eq:pde-sre},
\[
\left\langle
\nabla u(t,v),\nabla\phi(v)
\right\rangle=0,
\qquad
(t,v)\in[0,T)\times\partial\mathcal D .
\]
Moreover, the flatness condition \eqref{eq:flatness} implies that
the measure $\mathrm d\zeta_s^{t,v}$ is supported on
\[
\{s\in[t,T]:\bm V_s^{t,v}\in\partial\mathcal D\}.
\]
Therefore, the boundary integral vanishes.

By continuity, the first equation in \eqref{eq:pde-sre} also holds on
$[0,T)\times\partial\mathcal D$. Hence, we obtain
\[
\left(
\partial_\varsigma u+\mathcal L u
\right)
(\varsigma,\bm V_\varsigma^{t,v})
=
-\mathcal H
\left(
\varsigma,
\bm V_\varsigma^{t,v},
Z_\varsigma^{t,v}
\right).
\]
Since the terminal condition gives
\[
Y_T^{t,v}
=
u(T,\bm V_T^{t,v})
=
0,
\]
we conclude that
\[
Y_s^{t,v}
=
\int_s^T
\mathcal H
\left(
\varsigma,
\bm V_\varsigma^{t,v},
Z_\varsigma^{t,v}
\right)
\,\mathrm d\varsigma
-
\int_s^T
\left\langle
Z_\varsigma^{t,v},
\mathrm dW_\varsigma
\right\rangle .
\]
Hence, $(Y^{t,v},Z^{t,v})$ is an adapted solution of the quadratic BSDE
\eqref{eq:flow_quadratic}.
\end{proof}

Before constructing the viscosity solution of \eqref{eq:pde-sre} from the quadratic BSDE \eqref{eq:flow_quadratic}, we first recall the definition of viscosity solutions for semilinear parabolic equations with Neumann boundary conditions.

\begin{definition}
For a test function
\(\psi\in C^{1,2}([0,T]\times\overline{\mathcal D})\), define the interior and boundary operators by
\[
\begin{aligned}
\mathfrak R[\psi,u](t,v)
&:=
-\partial_t\psi(t,v)
-\mathcal L\psi(t,v)
-\mathcal H
\bigl(
t,v,\bm G(v)^*\nabla\psi(t,v)
\bigr),\\
\mathfrak N[\psi](t,v)
&:=
-\langle
\nabla\psi(t,v),
\nabla\phi(v)
\rangle .
\end{aligned}
\]

An upper semicontinuous (USC) function
\(u:[0,T]\times\overline{\mathcal D}\to\mathbb R\)
is called a viscosity subsolution of \eqref{eq:pde-sre} if
\(u(T,v)\le0\) for all \(v\in\overline{\mathcal D}\), and for every
\(\psi\in C^{1,2}([0,T]\times\overline{\mathcal D})\)
such that \(u-\psi\) attains a local maximum at
\((t_0,v_0)\in[0,T)\times\overline{\mathcal D}\), the following conditions hold:
\begin{itemize}
    \item If \(v_0\in\mathcal D\), then
    \[
    \mathfrak R[\psi,u](t_0,v_0)\le0.
    \]
    \item If \(v_0\in\partial\mathcal D\), then
    \[
    \min
    \left\{
    \mathfrak R[\psi,u](t_0,v_0),
    \mathfrak N[\psi](t_0,v_0)
    \right\}
    \le0.
    \]
\end{itemize}
A lower semicontinuous (LSC) function
$u:[0,T]\times\overline{\mathcal D}\to\mathbb R$
is called a viscosity supersolution of \eqref{eq:pde-sre}
if $u(T,v)\geq0$ for all $v\in\overline{\mathcal D}$, and for every
$\psi\in C^{1,2}([0,T]\times\overline{\mathcal D})$ such that
$u-\psi$ attains a local minimum at
$(t_0,v_0)\in[0,T)\times\overline{\mathcal D}$, the following conditions hold:
\begin{itemize}
    \item If $v_0\in\mathcal D$, then
    \[
        \mathfrak R[\psi,u](t_0,v_0)\geq0;
    \]
    \item If $v_0\in\partial\mathcal D$, then
    \[
        \max\left\{
            \mathfrak R[\psi,u](t_0,v_0),
            \mathfrak N[\psi](t_0,v_0)
        \right\}
        \geq0.
    \]
\end{itemize}

A continuous function
$u:[0,T]\times\overline{\mathcal D}\to\mathbb R$
is called a viscosity solution of \eqref{eq:pde-sre}
if it is both a viscosity subsolution and a viscosity supersolution.
\end{definition}

\begin{theorem}
\label{thm:u-continuity}
For each $(t,v)\in[0,T]\times\overline{\mathcal D}$, let
$(Y_s^{t,v},Z_s^{t,v})_{s\in[t,T]}$
be the adapted solution of the quadratic BSDE
\eqref{eq:flow_quadratic}. Define
\begin{equation} \label{eq:def_u_continuity}
    u(t,v):=Y_t^{t,v}.
\end{equation}
Then
\[
u\in C([0,T]\times\overline{\mathcal D}).
\]
\end{theorem}

\begin{proof}
By the standard Markovian property of the reflected FBSDE,
$Y_t^{t,v}$ is deterministic. Let
\[
(t_k,v_k)\rightarrow(t_0,v_0)
\]
in $[0,T]\times\overline{\mathcal D}$.

By the stability property of reflected SDEs
(see, for example, Corollary~4.56 in Pardoux and R{\u{a}}{\c{s}}canu \cite{pardoux2014stochastic}),
the reflected diffusion flows satisfy
\[
\sup_{s\in[t_k\vee t_0,T]}
|\bm V_s^{t_k,v_k}-\bm V_s^{t_0,v_0}|
\rightarrow0
\]
in probability.

To compare the BSDEs defined on different time intervals, we extend them
to the common interval $[0,T]$. Define
\[
\overline{\bm V}_s^k
=
\begin{cases}
v_k, &0\le s<t_k,\\
\bm V_s^{t_k,v_k},&t_k\le s\le T,
\end{cases}
\]
and
\[
\overline f_k(s,z)
=
\mathbf 1_{\{s\ge t_k\}}
\mathcal H(s,\overline{\bm V}_s^k,z).
\]

By Proposition~\ref{prop:cone-generator}, the generators
$\overline f_k$ satisfy uniform quadratic-growth and continuity estimates in $z$. Moreover,
\[
\overline f_k(\cdot,z)
\rightarrow
\overline f_0(\cdot,z)
\]
in
$\mathrm ds\otimes\mathrm d\mathbb P$-measure for every $z$.

Let
$(\overline Y^k,\overline Z^k)$ denote the solution of
\[
\overline Y_s^k
=
\int_s^T
\overline f_k(\varsigma,\overline Z_\varsigma^k)
\,\mathrm d\varsigma
-
\int_s^T
\langle\overline Z_\varsigma^k,\mathrm dW_\varsigma\rangle .
\]
By uniqueness,
\[
\overline Y_{t_k}^k
=
Y_{t_k}^{t_k,v_k}
=
u(t_k,v_k).
\]

The stability result for quadratic BSDEs (see, for example, Lemma 7.3.4 in Zhang \cite{zhang2017backward}) yields
\[
\lim_{k\to\infty}
\mathbb E
\left[
\sup_{0\le s\le T}
|\overline Y_s^k-\overline Y_s^0|^p
\right]
=0, \qquad  p \geq 1.
\]

Therefore,
\[
\begin{aligned}
|u(t_k,v_k)-u(t_0,v_0)|
&=
|\overline Y_{t_k}^k-\overline Y_{t_0}^0|
\le
\sup_{0\le s\le T}
|\overline Y_s^k-\overline Y_s^0|
+
|\overline Y_{t_k}^0-\overline Y_{t_0}^0|.
\end{aligned}
\]

The first term converges to zero in $L^p$, while the second term converges
to zero due to the continuity of the sample paths of $\overline Y^0$ and
$t_k\to t_0$. Hence,
\[
u(t_k,v_k)\rightarrow u(t_0,v_0),
\]
which proves the continuity of $u$.
\end{proof}

\begin{theorem}[Existence of a viscosity solution]
\label{thm:sre-to-viscosity}
The function $u$ defined by \eqref{eq:def_u_continuity} is a viscosity
solution of the parabolic terminal-boundary value problem
\eqref{eq:pde-sre}.
\end{theorem}

\begin{proof}
By Theorem~\ref{thm:u-continuity},
$u\in C([0,T]\times\overline{\mathcal D})$.
Moreover, since the terminal condition of
\eqref{eq:flow_quadratic} is zero,
\[
    u(T,v)=Y_T^{T,v}=0,
    \qquad v\in\overline{\mathcal D}.
\]

We first recall the flow property. By uniqueness of the quadratic BSDE and the Markov property of the reflected diffusion,
\begin{equation}\label{eq:dpp_flow}
Y_s^{t,v}
=
Y_s^{s,\bm V_s^{t,v}}
=
u(s,\bm V_s^{t,v}),
\qquad s\in[t,T],
\quad\text{a.s.}
\end{equation}
By continuity, this identity holds simultaneously
for all $s\in[t,T]$, almost surely.
Consequently, for any stopping time $\tau$ taking
values in $[t,T]$,
\begin{equation}\label{eq:local_bsde}
\begin{aligned}
    Y_s^{t,v}
    ={}&u(\tau,\bm V_\tau^{t,v})+\int_s^\tau
    \mathcal H(\varsigma,\bm V_\varsigma^{t,v},
    Z_\varsigma^{t,v})\,\mathrm d\varsigma-\int_s^\tau
    \langle Z_\varsigma^{t,v},\mathrm dW_\varsigma\rangle,
    \qquad t\leq s\leq\tau.
\end{aligned}
\end{equation}

We first verify the viscosity subsolution property.
Let $(t_0,v_0)\in[0,T)\times\overline{\mathcal D}$
and let $\psi\in C^{1,2}([0,T]\times\overline{\mathcal D})$
be such that $u-\psi$ attains a local maximum at
$(t_0,v_0)$. Without loss of generality, assume that
\[
    u(t_0,v_0)=\psi(t_0,v_0).
\]
Choose a sufficiently small relative open neighborhood $U$
of $(t_0,v_0)$ with $\overline U\subset[0,T)\times\overline{\mathcal D}$
such that
\begin{equation}\label{eq:touch_sub}
    u(s,v)\leq\psi(s,v),
    \qquad (s,v)\in\overline U.
\end{equation}

Suppose, by contradiction, that
\[
    \mathfrak R[\psi,u](t_0,v_0)>0
\]
if $v_0\in\mathcal D$, or
\[
    \min\{
        \mathfrak R[\psi,u](t_0,v_0),
        \mathfrak N[\psi](t_0,v_0)
    \}>0
\]
if $v_0\in\partial\mathcal D$.

By continuity, after shrinking $U$ if necessary,
there exists $\varepsilon_\psi>0$ such that
\begin{equation}\label{eq:sub_interior_strict}
\begin{aligned}
    -(\partial_s\psi+\mathcal L\psi)(s,v)
    \geq{}&
    \mathcal H(s,v,\bm G(v)^*\nabla\psi(s,v))
    +\varepsilon_\psi,
    \qquad (s,v)\in\overline U.
\end{aligned}
\end{equation}
In the interior case, we further choose $U$ so that
its closure does not intersect
$[0,T]\times\partial\mathcal D$.
In the boundary case, we also have
\begin{equation}\label{eq:sub_boundary_N}
    -\langle\nabla\psi(s,v),\nabla\phi(v)\rangle
    \geq\varepsilon_\psi,
    \qquad
    (s,v)\in\overline U\cap
    ([0,T)\times\partial\mathcal D).
\end{equation}

For notational simplicity, write
\[
    (\bm V_s,\zeta_s,Y_s,Z_s)
    :=
    (\bm V_s^{t_0,v_0},\zeta_s^{t_0,v_0},
    Y_s^{t_0,v_0},Z_s^{t_0,v_0}),
\]
and define the stopping time
\[
    \tau:=
    \inf\left\{
        s\in[t_0,T]:(s,\bm V_s)\notin U
    \right\}\wedge T.
\]
By continuity, $\tau>t_0$ almost surely and
$(\tau,\bm V_\tau)\in\overline U$.

Set
\[
    \widehat Y_s:=\psi(s,\bm V_s),
    \qquad
    \widehat Z_s:=
    \bm G(\bm V_s)^*\nabla\psi(s,\bm V_s).
\]
Applying It\^o's formula yields
\begin{equation} \label{quadratic-ito}
    \begin{aligned}
            \widehat Y_s
    ={}&\psi(\tau,\bm V_\tau)
    +\int_s^\tau\widehat f_\varsigma^\psi
        \,\mathrm d\varsigma+\int_s^\tau\widehat g_\varsigma^\psi
        \,\mathrm d\zeta_\varsigma
    -\int_s^\tau
        \langle\widehat Z_\varsigma,
        \mathrm dW_\varsigma\rangle,
    \end{aligned}
\end{equation}
where
\[
    \widehat f_s^\psi
    :=
    -(\partial_s\psi+\mathcal L\psi)(s,\bm V_s),
    \qquad
    \widehat g_s^\psi
    :=
    -\langle\nabla\psi(s,\bm V_s),
        \nabla\phi(\bm V_s)\rangle.
\]
By \eqref{eq:sub_interior_strict},
\[
    \widehat f_s^\psi
    -
    \mathcal H(s,\bm V_s,\widehat Z_s)
    \geq\varepsilon_\psi,
    \qquad s\in[t_0,\tau].
\]
Moreover, the reflection term vanishes in the
interior case, while in the boundary case,
\eqref{eq:sub_boundary_N} and the support property
of $\mathrm d\zeta$ imply that
$\widehat g_s^\psi\geq\varepsilon_\psi$,
$\mathrm d\zeta_s$-a.e. Hence, in both cases,
\[
    \int_s^\tau
        \widehat g_\varsigma^\psi
        \,\mathrm d\zeta_\varsigma
    \geq0.
\]

Define
\[
    \chi_s^\psi
    :=
    \begin{cases}
    \displaystyle
    \frac{
        \mathcal H(s,\bm V_s,\widehat Z_s)
        -
        \mathcal H(s,\bm V_s,Z_s)
    }{|\widehat Z_s-Z_s|^2}
    (\widehat Z_s-Z_s),
    &\qquad  \widehat Z_s\neq Z_s,\\[8pt]
    0, &\qquad  \widehat Z_s=Z_s.
    \end{cases}
\]
By the weighted Lipschitz estimate in Proposition~\ref{prop:cone-generator},
\[
    |\chi_s^\psi|
    \leq C(1+|\widehat Z_s|+|Z_s|).
\]
Since $Y$ is bounded and $\mathcal H$ has quadratic
growth, the stochastic integral of $Z$ with respect
to $W$ is a BMO martingale. Together with the
boundedness of $\widehat Z$, this implies that
\[
    \int_{t_0}^{\cdot\wedge\tau}
        \langle\chi_s^\psi,\mathrm dW_s\rangle
\]
is also a BMO martingale.

We may therefore define an equivalent probability
measure $\mathbb P^\psi$ by
\[
    \frac{\mathrm d\mathbb P^\psi}{\mathrm d\mathbb P}
    :=
    \mathcal E\left(
        \int_{t_0}^{\cdot\wedge\tau}
        \langle\chi_s^\psi,\mathrm dW_s\rangle
    \right)_T.
\]
By Girsanov's theorem,
\[
    W_s^\psi
    :=
    W_s-\int_{t_0}^{s\wedge\tau}
        \chi_\varsigma^\psi\,\mathrm d\varsigma,
    \qquad s\in[t_0,T],
\]
is a Brownian motion under $\mathbb P^\psi$ on
$[t_0,T]$.

Subtracting \eqref{eq:local_bsde} from the BSDE \eqref{quadratic-ito} satisfied by $\widehat Y$ yields
\begin{align*}
    \widehat Y_s-Y_s
    ={}&\widehat Y_\tau-Y_\tau+\int_s^\tau
        \bigl[
            \widehat f_\varsigma^\psi
            -
            \mathcal H(\varsigma,\bm V_\varsigma,
                \widehat Z_\varsigma)
        \bigr]\,\mathrm d\varsigma\\
    &+\int_s^\tau
        \widehat g_\varsigma^\psi\,\mathrm d\zeta_\varsigma - \int_s^\tau
        \langle\widehat Z_\varsigma-Z_\varsigma,
        \mathrm dW_\varsigma^\psi\rangle.
\end{align*}

Let $(\tau_j)_{j\geq1}$ be a localizing sequence
for the stochastic integral under $\mathbb P^\psi$,
with $t_0\leq\tau_j\leq\tau$ and
$\tau_j\uparrow\tau$ almost surely.
Using the preceding inequalities and taking
conditional expectations, we obtain
\[
    \widehat Y_{t_0}-Y_{t_0}
    \geq
    \mathbb E^{\mathbb P^\psi}\left[
        \widehat Y_{\tau_j}-Y_{\tau_j}
        +\varepsilon_\psi(\tau_j-t_0)
        \,\middle|\,\mathcal F_{t_0}
    \right].
\]
By boundedness and continuity of $Y$ and
$\widehat Y$, letting $j\to\infty$ gives
\[
    \widehat Y_{t_0}-Y_{t_0}
    \geq
    \mathbb E^{\mathbb P^\psi}\left[
        \widehat Y_\tau-Y_\tau
        +\varepsilon_\psi(\tau-t_0)
        \,\middle|\,\mathcal F_{t_0}
    \right].
\]
By \eqref{eq:touch_sub} and \eqref{eq:dpp_flow},
\[
    \widehat Y_\tau-Y_\tau
    =
    \psi(\tau,\bm V_\tau)-u(\tau,\bm V_\tau)
    \geq0.
\]
Since $\tau>t_0$ almost surely, it follows that
\[
    \widehat Y_{t_0}-Y_{t_0}
    \geq
    \varepsilon_\psi
    \mathbb E^{\mathbb P^\psi}\left[
        \tau-t_0
        \,\middle|\,\mathcal F_{t_0}
    \right]
    >0.
\]
Consequently,
\[
    \psi(t_0,v_0)>u(t_0,v_0),
\]
contradicting the touching condition.
Thus $u$ is a viscosity subsolution.

The viscosity supersolution property follows
analogously by considering a test function touching
$u$ from below and reversing the above inequalities.
Therefore, $u$ is a viscosity solution of
\eqref{eq:pde-sre}.
\end{proof}

We next introduce a change of variables inspired by
Kobylanski \cite{kobylanski2000backward}.
Since the original Hamiltonian is independent of the
solution variable, we transform the equation into an
equivalent Neumann problem whose Hamiltonian satisfies
a strict monotonicity condition required for the
comparison argument in Theorem \ref{thm:comparison-transformed-neumann}. 

\begin{theorem}\label{thm:kobylanski-transform-neumann}
Assume that $\mathcal H$ satisfies the quadratic-growth
and weighted Lipschitz estimates of
Proposition~\ref{prop:cone-generator}.
Let $J\subset\mathbb R$ be a bounded interval and fix
$M>0$ such that $J\subset[-M,M]$.
Then there exist constants $K>0$ and $\kappa>0$,
and a strictly increasing $C^\infty$-diffeomorphism
$\Psi:\mathbb R\to\mathbb R$ with $\Psi(0)=0$,
depending only on $M$, $T$, and the structural constants
in Proposition~\ref{prop:cone-generator},
such that, for any viscosity subsolution, supersolution,
or solution $u$ of \eqref{eq:pde-sre} taking values in $J$,
the transformation
\[
    u(t,v)=e^{K(T-t)}\Psi(\widetilde u(t,v))
\]
maps $u$ into a viscosity subsolution, supersolution,
or solution, respectively, of the transformed
Neumann problem
\begin{equation} \label{transformed-pde}
    \left\{
\begin{aligned}
\partial_t\widetilde u+\mathcal L\widetilde u
+\widetilde{\mathcal H}_K
\bigl(t,v,\widetilde u,\bm G(v)^*\nabla\widetilde u\bigr)&=0,
&& (t,v)\in[0,T)\times\mathcal D,\\
\langle\nabla\widetilde u(t,v),\nabla\phi(v)\rangle&=0,
&& (t,v)\in[0,T)\times\partial\mathcal D,\\
\widetilde u(T,v)&=0,
&& v\in\overline{\mathcal D}.
\end{aligned}
\right.
\end{equation}
Let
\[
    I_M:=\Psi^{-1}([-M,M]).
\]
Then $\widetilde u$ takes values in $I_M$, and
the transformed Hamiltonian satisfies
\[
\widetilde{\mathcal H}_K(t,v,\xi_1,z)
-
\widetilde{\mathcal H}_K(t,v,\xi_2,z)
\leq
-\kappa(\xi_1-\xi_2)(1+|z|^2)
\]
for all
$(t,v,z)\in[0,T]\times\overline{\mathcal D}\times\mathbb R^n$
and $\xi_1,\xi_2\in I_M$ with $\xi_1\geq\xi_2$.
Moreover, after $\Psi$ and $K$ are fixed,
$\widetilde{\mathcal H}_K$ satisfies quadratic-growth
and weighted Lipschitz estimates analogous to those in
Proposition~\ref{prop:cone-generator}.
\end{theorem}

The proof is deferred to Appendix~\ref{app-kobi}.

\begin{remark}
In contrast to the Cauchy problem on the whole space
considered by Kobylanski \cite{kobylanski2000backward},
we study a Neumann problem on a bounded domain.
The requirement that $J$ contain the range of $u$
introduces no additional restriction in the
bounded-solution framework. Moreover, choosing a common
interval $J$ for two bounded viscosity sub- and
supersolutions allows us to apply the same transformation
to both.

Unlike the transformation argument in Kobylanski \cite{kobylanski2000backward}, which relies on differentiability assumptions on the Hamiltonian, our argument derives the required strict monotonicity from the quadratic-growth and weighted Lipschitz estimates in Proposition~\ref{prop:cone-generator}, using the almost-everywhere differentiability of a one-dimensional Lipschitz composition. In particular, differentiability of the Hamiltonian with respect to
the gradient variable is not required.
\end{remark}

The uniqueness of bounded continuous viscosity solutions to
\eqref{eq:pde-sre} follows from the comparison principle below. Its proof is based on a doubling-of-variables argument adapted to the homogeneous Neumann boundary condition and given in Appendix \ref{app-comparison}.

\begin{theorem}[Comparison principle]
\label{thm:comparison-transformed-neumann}
Let
\[
\widetilde u\in
\mathrm{USC}([0,T]\times\overline{\mathcal D}),
\qquad
\widetilde w\in
\mathrm{LSC}([0,T]\times\overline{\mathcal D}),
\]
be bounded viscosity sub- and supersolutions of the transformed Neumann problem \eqref{transformed-pde}. Assume that the terminal conditions are ordered, namely $\widetilde u(T,v)\leq \widetilde w(T,v)$ for all $v\in\overline{\mathcal D}$ and the ranges of $\widetilde u$ and $\widetilde w$ are contained in a bounded interval $I_M$. Suppose that, for some $\kappa>0$,
\[
\widetilde{\mathcal H}_K(t,v,\xi_1,z)
-
\widetilde{\mathcal H}_K(t,v,\xi_2,z)
\leq
-\kappa(\xi_1-\xi_2)(1+|z|^2),
\]
for all $(t,v,z)\in[0,T]\times\overline{\mathcal D}
\times\mathbb R^n$ and $\xi_1,\xi_2\in I_M$
with $\xi_1\geq\xi_2$, and that
$\widetilde{\mathcal H}_K$ satisfies the weighted continuity estimate
\[
\begin{aligned}
&\quad
\left|
\widetilde{\mathcal H}_K(t_1,v_1,\xi,z_1)
-
\widetilde{\mathcal H}_K(t_2,v_2,\xi,z_2)
\right|
\\
&\leq
\omega_f(|t_1-t_2|+|v_1-v_2|)
\bigl(1+|z_1|^2+|z_2|^2\bigr)
+
C(1+|z_1|+|z_2|)|z_1-z_2|,
\end{aligned}
\]
for all $(t_i,v_i,z_i)\in[0,T]\times
\overline{\mathcal D}\times\mathbb R^n$, $i=1,2$,
and $\xi\in I_M$, where $\omega_f$ is a modulus
of continuity. Then
\[
\widetilde u(t,v)\leq\widetilde w(t,v),
\qquad
(t,v)\in[0,T]\times\overline{\mathcal D}.
\]
\end{theorem}

Consequently, the transformed Neumann problem
\eqref{transformed-pde} admits at most one bounded
continuous viscosity solution taking values in $I_M$.
For any two bounded continuous viscosity solutions of
\eqref{eq:pde-sre}, we may choose a common interval $J$
containing their ranges and apply the same transformation
in Theorem~\ref{thm:kobylanski-transform-neumann}.
The comparison principle, together with the strict
monotonicity of $\Psi$, then yields their uniqueness.

\section{A Deep BSDE Framework for Markovian BSDEs with Reflected Factors}

In this section, we develop a deep BSDE method for the logarithmically transformed equations \eqref{log-positive} and \eqref{log-negative} driven by reflected Markovian factors, building on the method of Han, Jentzen, and E \cite{han2018solving}. Related deep-learning approaches for Neumann boundary problems have employed reflected stochastic processes and FBSDE-based neural networks; see, for example, \cite{boussange2023deep,deng2026novel,wang2026deep}. Our method focuses on stochastic Riccati equations
with reflected factors and combines reflected-path
simulation, geometry-aware boundary treatment,
neural approximation of the martingale integrand,
and efficient evaluation of the constrained
Hamiltonian.

\subsection{Reflected Ornstein--Uhlenbeck Factor Process}

Reflected SDEs can be simulated using Euler-type discretization schemes; see, for example, Gobet~\cite{gobet2001euler}. In our numerical implementation, we specialize the general
reflected factor dynamics \eqref{eq:RSDE} to an
$n$-dimensional reflected Ornstein--Uhlenbeck (ROU) process
on the closed ball $\overline{\mathcal D}$, where
\[
\mathcal D=\{v\in\mathbb R^n:|v|<R\},
\qquad R>0.
\]
This provides a simple and computationally tractable realization of the reflected-factor framework while preserving its essential features: the factor remains in the admissible domain and the reflection acts only at the boundary in the normal direction. Specifically, the factor process satisfies
\[
\mathrm{d}V_t
=
\beta_0(\alpha_0-V_t)\,\mathrm{d}t
+
\sigma_0\,\mathrm{d}W_t
-
\gamma(V_t)\,\mathrm{d}L_t,
\qquad
V_0=v_0\in\overline{\mathcal D},
\]
where $\beta_0\in\mathbb R^{n\times n}$ is the mean-reversion matrix, $\alpha_0\in\mathbb R^n$ is the long-run mean,
$\sigma_0\in\mathbb R^{n\times n}$ is the diffusion matrix, and
$W$ is an $n$-dimensional Brownian motion. Here $L$ is the nondecreasing boundary regulator corresponding to $\zeta$ in \eqref{eq:RSDE}. Let
\[
\gamma(v)=\frac{v}{|v|},\qquad v\in\partial\mathcal D,
\]
denote the outward unit normal to $\partial\mathcal D$ at $v$.
The term $-\gamma(V_t)\,\mathrm{d}L_t$ thus realizes the
inward normal reflection that keeps $V_t$ in
$\overline{\mathcal D}$.
This formulation is consistent with \eqref{eq:RSDE}.
Indeed, setting
\[
\phi(v)=\frac{R(R^2-|v|^2)}{R^2+|v|^2},
\]
we have $\mathcal D=\{\phi>0\}$,
$\partial\mathcal D=\{\phi=0\}$, and
$\nabla\phi(v)=-\gamma(v)$ for $v\in\partial\mathcal D$.

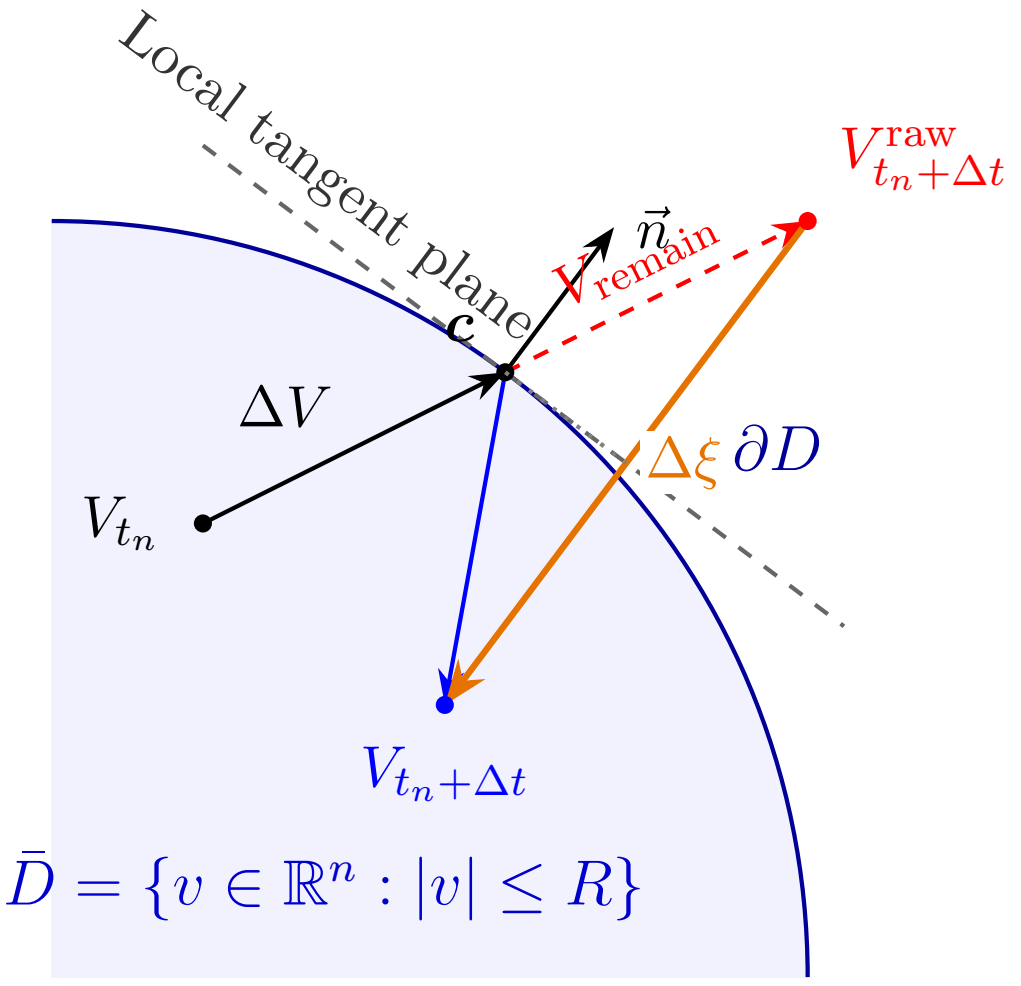
\begin{figure}[htbp]
\centering
\begin{tikzpicture}[scale=1.0, >=Stealth]
    
    % 1. Setup the Domain D and its boundary
    % Use a clipping scope to only color the inside of the domain
    \begin{scope}
        \clip (0, 0) rectangle (6, 6);
        \fill[blue!5] (0,0) circle (5);
    \end{scope}
    
    % Draw the boundary arc (Radius = 5)
    \draw[thick, blue!60!black] (5,0) arc (0:90:5);
    \node[text=blue!80!black, font=\large] at (1.8, 0.6) {$\bar{D} = \{ v \in \mathbb{R}^n : |v| \le R \}$};
    \node[text=blue!60!black, font=\large] at (4.8, 3.5) {$\partial D$};

    % 2. Exact Coordinate Calculations based on the algorithm
    % Center O=(0,0), Radius=5. Collision point c=(3,4) lies exactly on the boundary.
    \coordinate (c) at (3, 4);
    
    % Starting point Vn=(1,3). The ray travels (2,1) to hit c.
    \coordinate (Vn) at (1, 3);
    
    % V_raw = c + V_remain = (3,4) + (2,1) = (5,5)
    \coordinate (Vraw) at (5, 5);
    
    % Normal vector at c is (3,4), normalized to n=(0.6, 0.8)
    \coordinate (N) at ($ (c) + 1.2*(0.6, 0.8) $);
    
    % Reflection geometry calculation:
    % V_remain = (2,1). Inner product with n = 2(0.6) + 1(0.8) = 2.0
    % Reflected V_remain = V_remain - 2(2.0)n = (2,1) - 4(0.6, 0.8) = (2,1) - (2.4, 3.2) = (-0.4, -2.2)
    % V_Normal = c + Reflected = (3,4) + (-0.4, -2.2) = (2.6, 1.8)
    \coordinate (Vspec) at (2.6, 1.8);

    % Tangent line at c (direction is orthogonal to n, i.e., (-0.8, 0.6))
    \coordinate (T1) at ($ (c) + 2.5*(-0.8, 0.6) $);
    \coordinate (T2) at ($ (c) - 2.8*(-0.8, 0.6) $);

    % 3. Draw Tangent Plane and Normal Vector
    \draw[dashed, thick, black!60] (T1) -- (T2) 
        node[pos=0.15, above, sloped, text=black!80] {Local tangent plane};
    \draw[->, thick] (c) -- (N) node[right] {$\vec n$};

    % 4. Draw Trajectories
    % Step 1: Vn to c (Trajectory inside domain)
    \draw[->, thick] (Vn) -- (c) node[midway, above left] {$\Delta V$};
    
    % Step 2: c to Vraw (Overshooting segment)
    \draw[->, thick, dashed, red] (c) -- (Vraw) node[midway, sloped, above] {$V_{\text{remain}}$};
    
    % Step 3: c to Vspec (Reflected segment)
    \draw[->, thick, blue] (c) -- (Vspec);

    % 5. Draw Skorokhod Displacement (Correction vector)
    % From Vraw to Vspec. Note it is perfectly parallel to -n(c).
    \draw[->, very thick, orange!90!black] (Vraw) -- (Vspec) 
        node[midway, right=2pt, fill=white, inner sep=1pt] {$\Delta \xi$};

    % 6. Draw Points and Text Labels
    \filldraw (Vn) circle (1.5pt) node[left=4pt] {$V_{t_n}$};
    \filldraw (c) circle (1.5pt) node[above left=2pt] {$\bm c$};
    \filldraw[red] (Vraw) circle (1.5pt) node[above right=2pt, text=red] {$V_{t_n+\Delta t}^{\rm raw}$};
    \filldraw[blue] (Vspec) circle (1.5pt) node[below=4pt, text=blue] {$V_{t_n+\Delta t}$};

    % Optional: Geometric right-angle indicator for reflection symmetry
    % The line Vraw->Vspec intersects the tangent plane exactly at their midpoint.
    \coordinate (Mid) at ($ (Vraw)!0.5!(Vspec) $);
    \draw[dotted, thick, gray] (c) -- (Mid);

\end{tikzpicture}
\caption{Illustration of the discrete-time normal reflection scheme.}
\label{fig:Normal_reflection}
\end{figure}

For numerical simulation, we first perform an Euler--Maruyama step,
\[
V_{t+\Delta t}^{\rm raw}
=
V_t
+
\beta_0(\alpha_0-V_t)\Delta t
+
\sigma_0\Delta W_t.
\]
If $|V_{t+\Delta t}^{\rm raw}|\leq R$, the candidate is accepted. Otherwise, we locate the point $c$ at which the
segment joining $V_t$ and $V_{t+\Delta t}^{\rm raw}$
exits the closed ball $\overline{\mathcal D}$. Let
\[
\vec n=\frac{\bm c}{|\bm c|}
\]
be the corresponding outward normal vector. The overshooting segment is reflected across the tangent plane:
\[
V_{t+\Delta t}
=
\bm c
+
\bigl(V_{t+\Delta t}^{\rm raw}-\bm c\bigr)
-
2\bigl(
(V_{t+\Delta t}^{\rm raw}-\bm c)\cdot\vec n
\bigr)\vec n .
\]
A final radial projection is applied if necessary to ensure
$V_{t+\Delta t}\in\overline{\mathcal D}$. This reflection-correction scheme is used throughout our numerical experiments. Its geometric construction is illustrated in Figure~\ref{fig:Normal_reflection}.

\subsection{Neural Network Approximation}

Let $0=t_0<t_1<\cdots<t_N=T$ be a time partition,
with $\Delta t_i=t_{i+1}-t_i$ and $\Delta W_{t_i}=W_{t_{i+1}}-W_{t_i}$. For notational simplicity, we write $(Y,Z)$ for
$(\check P_\pm,\check\Lambda_\pm)$ and denote the
corresponding generator by $\mathcal H$. We now construct a neural-network approximation for the quadratic BSDE associated with the reflected factor process:
\begin{equation*}\label{deep-bsde}
\left\{
\begin{aligned}
\mathrm dY_t
&=
-\mathcal H(t,V_t,Z_t)\,\mathrm dt
+
\langle Z_t,\mathrm dW_t\rangle,
\qquad t\in[0,T),\\
Y_T&=0,
\end{aligned}
\right.
\end{equation*}
where $\mathcal H$ satisfies the quadratic-growth and weighted Lipschitz
estimates of Proposition~\ref{prop:cone-generator}. At each time step
$t_i$, a neural network produces a raw approximation
$Z_{t_i}^{\rm raw}$ of the martingale integrand.

The reflected setting requires particular attention near the boundary. For the ROU model, when the Markovian decoupling field \(u\) is sufficiently smooth, the martingale component admits the representation
\[
Z_t=\sigma_0^*\nabla u(t,V_t).
\]
This identity motivates the following boundary-aware correction. Assuming that $\sigma_0$ is invertible, the homogeneous Neumann condition
is equivalently written as
\[
\left\langle
Z_t,\sigma_0^{-1}\vec n(V_t)
\right\rangle=0,
\qquad V_t\in\partial\mathcal D,
\]
where $\vec n$ denotes the outward unit normal. To incorporate this condition into the network architecture, we introduce a geometry-aware correction layer. For convenience, we extend the outward unit normal to $\overline{\mathcal D}$ by setting
\[
\vec n(v)=
\begin{cases}
v/|v|, & v\neq 0,\\
0, & v=0.
\end{cases}
\]
Let
\[
\omega_0(x)
=
\max\left\{
\left(1-\frac{x}{d_{\rm th}}\right)^3,0
\right\},
\]
where $d_{\rm th} \in (0,R)$ is a prescribed boundary threshold, and define
\[
\Delta Z_{t_i}
=
\omega_0\bigl(\operatorname{dist}(V_{t_i},\partial\mathcal D)\bigr)
\left\langle
Z_{t_i}^{\rm raw},
\sigma_0^{-1}\vec n(V_{t_i})
\right\rangle
\sigma_0^*\vec n(V_{t_i})
\]
followed by
\begin{equation*}\label{eq:projection}
Z_{t_i}
=
Z_{t_i}^{\rm raw}-\Delta Z_{t_i}.
\end{equation*}
The correction enforces
\[
\left\langle
Z_{t_i},\sigma_0^{-1}\vec n(V_{t_i})
\right\rangle=0,
\]
whenever $V_{t_i}\in\partial\mathcal D$. Away from the boundary, the
cutoff smoothly reduces the correction. This construction is related to
analytical correction layers for constrained neural-network outputs in Dalal et al.
\cite{dalal2018safe} and geometry-based treatments of boundary conditions in Sukumar
and Srivastava \cite{sukumar2022exact}.

Let the reflected Euler scheme for the forward factor be written as
\[
V_{t_{i+1}}
=
V_{t_i}
+
\beta_0(\alpha_0-V_{t_i})\Delta t_i
+
\sigma_0\Delta W_{t_i}
+
\Delta\xi_{t_i},
\]
where $\Delta\xi_{t_i}$ denotes the discrete reflection correction.
Given $Z_{t_i}$, we propagate the backward component by
\begin{equation*}\label{eq:compensated_EM}
\begin{aligned}
Y_{t_{i+1}}
={}&
Y_{t_i}
-
\mathcal H(t_i,V_{t_i},Z_{t_i})\Delta t_i
+
\langle Z_{t_i},\Delta W_{t_i}\rangle +
\left\langle
(\sigma_0^*)^{-1}Z_{t_i},
\Delta\xi_{t_i}
\right\rangle .
\end{aligned}
\end{equation*}
Although the reflection term in Itô's formula \eqref{eq:Ito}
vanishes in continuous time under the homogeneous
Neumann condition, its discrete counterpart need
not vanish due to the left-endpoint evaluation
of the gradient. We retain this term to account
for the displacement induced by the numerical
reflection.

The initial value $Y_0$ and the neural-network parameters are trained by minimizing
\begin{equation*}\label{eq:loss_function}
\operatorname{Loss}(Y_0,\theta)
=
\mathbb E\bigl[|Y_T|^2\bigr]
+
\lambda_{\rm wall}
\mathbb E
\left[
\sum_{i=0}^{N-1}
|\Delta Z_{t_i}|^2
\right],
\end{equation*}
where the first term enforces the terminal condition and the second encourages the raw network output to satisfy the Neumann condition near the boundary. The parameter $\lambda_{\rm wall}\geq0$ controls the strength of the boundary
regularization.

\begin{remark}
    Compared with existing deep BSDE methods for Neumann problems~\cite{boussange2023deep,deng2026novel,wang2026deep}, our approach combines reflected-path simulation with a geometry-aware correction of the martingale integrand and PINN-inspired boundary regularization.
\end{remark}

\subsection{Numerical Evaluation of the Infimum Operator via FISTA}

For the numerical implementation, we focus on
the no-short-selling constraint $\Gamma=\mathbb R_+^m$,
under which the pointwise minimization problem in \eqref{log-positive} and \eqref{log-negative} 
is a strongly convex quadratic program. To obtain an efficient GPU-compatible solver, we employ the fast iterative shrinkage-thresholding algorithm (FISTA) of Beck and Teboulle~\cite{MR2486527}.

At each time step $t_i$ and factor state $V_{t_i}$, define
\[
A_i:=\sigma(t_i,V_{t_i})\sigma(t_i,V_{t_i})^*,
\qquad
b_{\pm,i}
:=
\widetilde\mu(t_i,V_{t_i})
+
\sigma(t_i,V_{t_i})\check\Lambda_{\pm,t_i}.
\]
The inner optimization problem is
\begin{equation*}\label{eq:inner_QP}
\pi_{\pm,i}^*
=
\arg\min_{\pi\in\Gamma}
\mathcal Q_{\pm,i}(\pi),
\qquad
\mathcal Q_{\pm,i}(\pi)
:=
\langle A_i\pi,\pi\rangle
\pm2\langle b_{\pm,i},\pi\rangle .
\end{equation*}
% The vector $b_{\pm,i}$ combines the instantaneous excess return with the
% hedging component generated by the neural-network approximation of the
% martingale integrand.

In the no-short-selling case,
\[
\Gamma=\mathbb R_+^m,
\]
the proximal operator of the indicator function of $\Gamma$ reduces to the
componentwise projection
\[
\mathcal T_{\geq0}(x):=\max(x,0).
\]
We initialize the iteration by projecting the unconstrained minimizer,
\[
\pi^{(0)}
=
\mathcal T_{\geq0}
\left(
\mp A_i^{-1}b_{\pm,i}
\right),
\qquad
y^{(0)}=\pi^{(0)},
\qquad
a_0=1.
\]
The FISTA iterations are then
\[
\pi^{(k+1)}
=
\mathcal T_{\geq0}
\left(
y^{(k)}
-
\eta_i
\nabla\mathcal Q_{\pm,i}(y^{(k)})
\right),
\]
\[
a_{k+1}
=
\frac{1+\sqrt{1+4a_k^2}}{2},
\]
and
\[
y^{(k+1)}
=
\pi^{(k+1)}
+
\frac{a_k-1}{a_{k+1}}
\left(
\pi^{(k+1)}-\pi^{(k)}
\right),
\]
where
\[
\nabla\mathcal Q_{\pm,i}(\pi)
=
2\left(
A_i\pi\pm b_{\pm,i}
\right).
\]

We use the adaptive step size
\[
\eta_i
=
\frac{c_\eta}
{\operatorname{Tr}(A_i)+\epsilon},
\qquad
0<c_\eta\leq\frac12,
\]
with a small safeguard $\epsilon>0$. Since
\[
2\lambda_{\max}(A_i)
\leq
2\operatorname{Tr}(A_i),
\]
this choice provides a conservative approximation to the reciprocal Lipschitz constant of the gradient and is convenient for batched GPU implementation.

After a prescribed number of FISTA iterations, the
resulting approximate minimizer $\widehat\pi_{\pm,i}$
is used to evaluate the corresponding BSDE generator. All matrix--vector operations and projections are performed simultaneously across Monte Carlo paths.

\subsection{Numerical Results}

We now report numerical results for the transformed BSDEs
\eqref{log-positive} and \eqref{log-negative}. The market coefficients are driven by a reflected Ornstein--Uhlenbeck factor process, following the factor-based stochastic-volatility modeling paradigm of Feng, Fouque, and Kumar~\cite{MR2985169}. The drift and volatility coefficients retain their
Markovian factor structure, while the factor process is normally reflected inside $\overline{\mathcal D}$ using the reflection scheme described above. 

As an illustrative high-dimensional numerical example, we consider a market with $n=25$ factors and $m=12$ risky assets under the no-short-selling constraint $\Gamma=\mathbb R_+^m$, equipped with a synthetic parameter specification designed to test the proposed solver. Under this specification, the boundary reflection mechanism is actively engaged: approximately $70\%$ of the simulated paths experience at least one boundary reflection, with an average of about $5$ reflections per path. All model specifications, synthetic coefficient matrices, heatmaps, reflection statistics, and solver implementations are publicly available at
\begin{center}
\url{https://github.com/huangzhecheng1996/DeepBSDE_Reflected_OU/}
\end{center}
to facilitate reproducibility.

\paragraph{Training dynamics.}

\begin{figure}[htbp]
    \centering
    \begin{minipage}{0.49\textwidth}
        \centering
        \includegraphics[width=\linewidth]{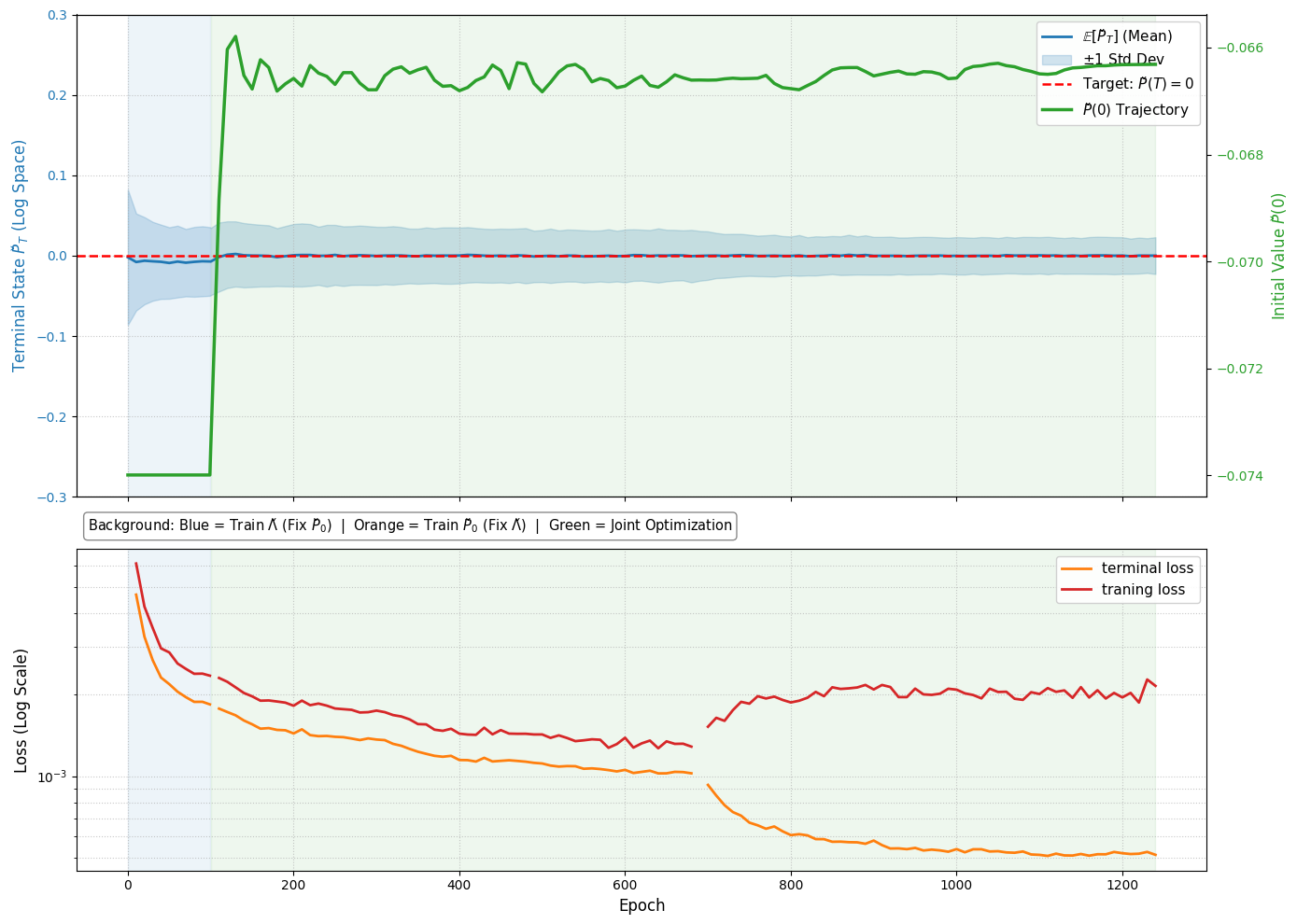}
        \caption*{$\check P_+$}
    \end{minipage}
    \hfill
    \begin{minipage}{0.49\textwidth}
        \centering
        \includegraphics[width=\linewidth]{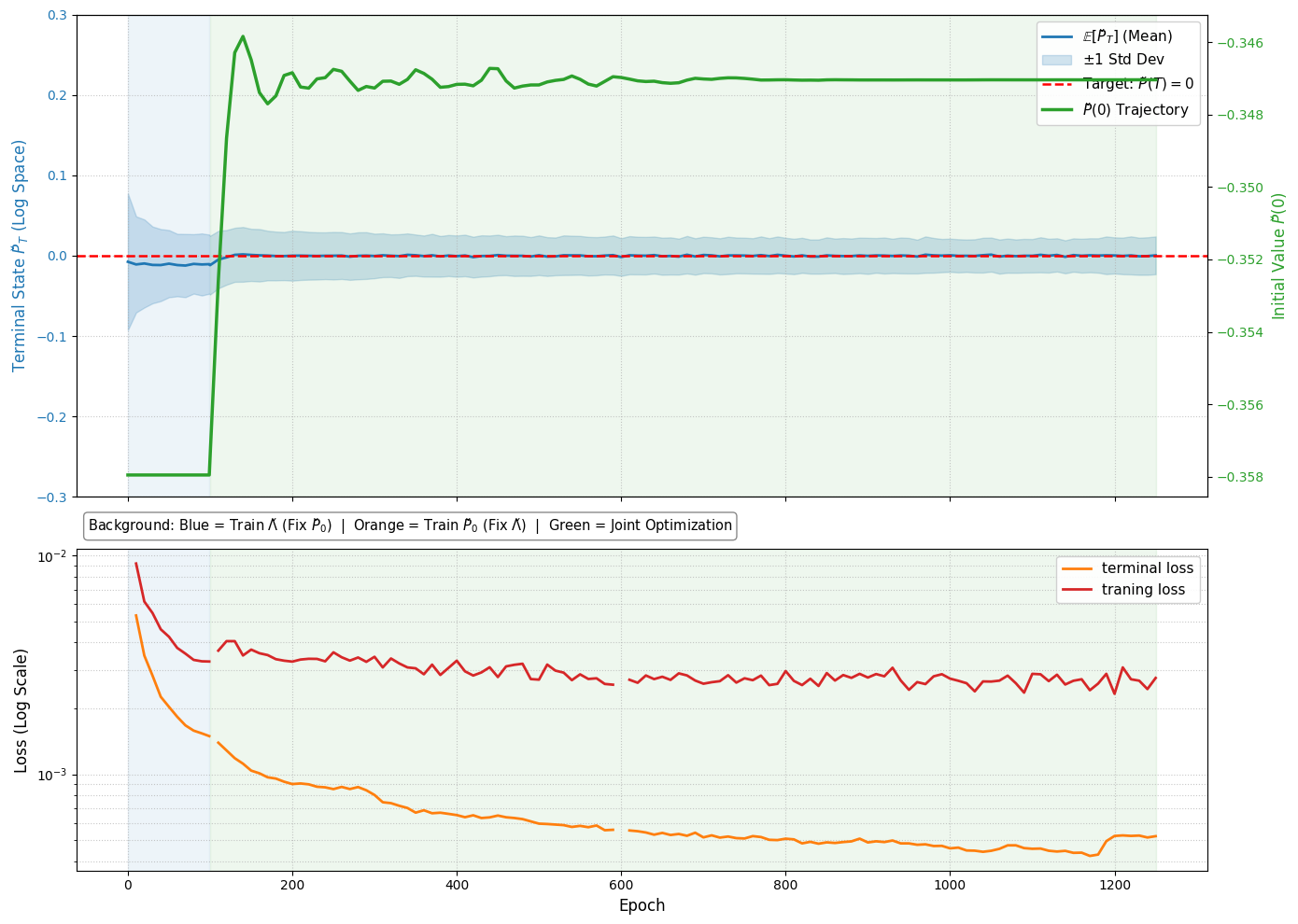}
        \caption*{$\check P_-$}
    \end{minipage}
    \caption{Training dynamics for the logarithmically transformed BSDEs.}
    \label{fig:training_history}
\end{figure}

Figure~\ref{fig:training_history} reports the training histories for the two logarithmically transformed equations, denoted by $\check P_+$ and $\check P_-$. The upper panels show the learned initial values together with the empirical terminal distributions, while the lower panels report the
training losses and terminal losses on a logarithmic scale. In both cases, the terminal distributions become increasingly concentrated around zero, indicating stable training behavior. The total training loss exhibits a moderate increase at a later stage, while the terminal loss continues to decrease. This behavior reflects the interplay
between terminal fitting and boundary regularization.

\paragraph{Out-of-sample evaluation.}

\begin{figure}[htbp]
    \centering
    \begin{minipage}{0.49\textwidth}
        \centering
        \includegraphics[width=\linewidth]{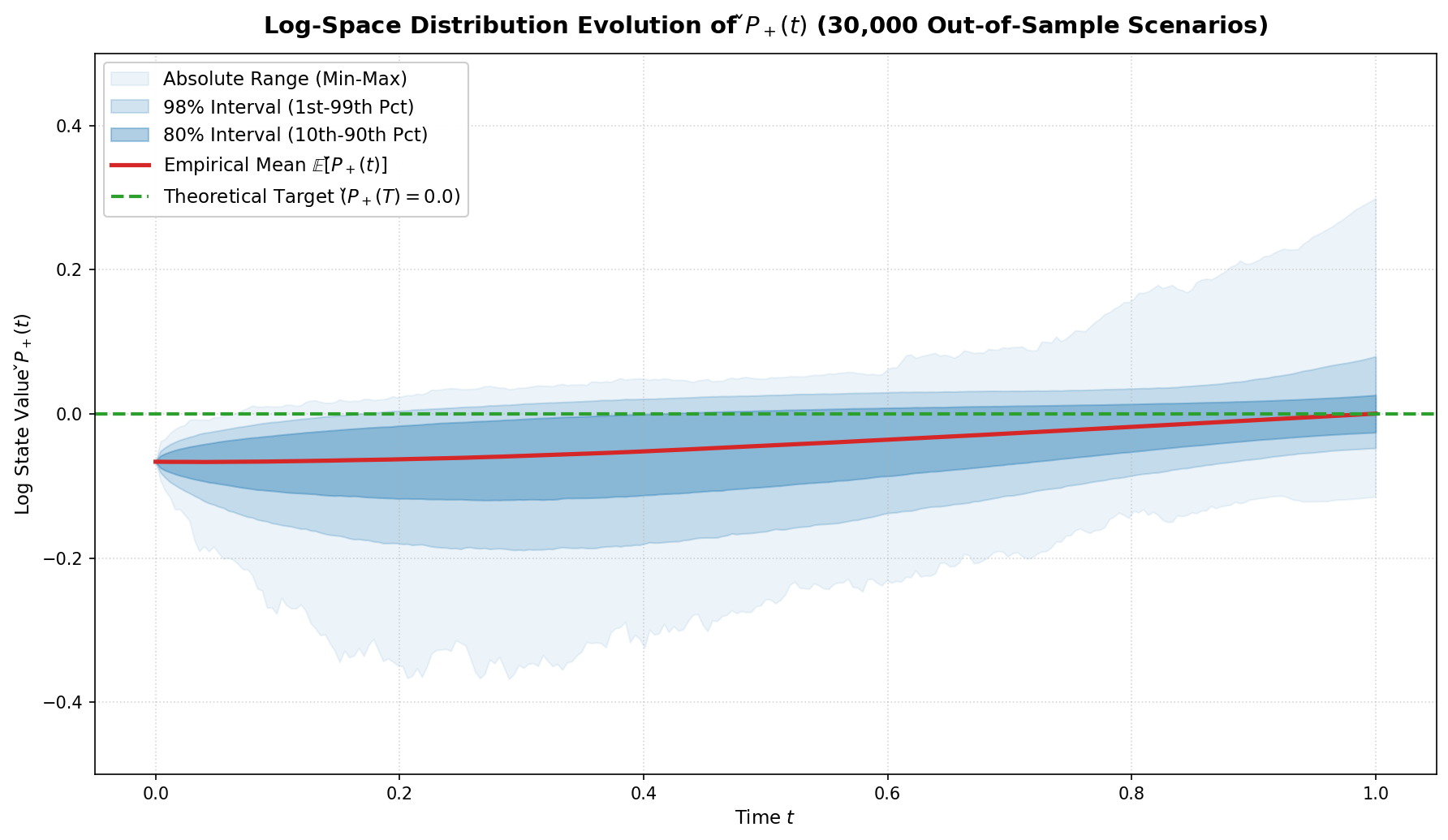}
        \caption*{$\check P_+$}
    \end{minipage}
    \hfill
    \begin{minipage}{0.49\textwidth}
        \centering
        \includegraphics[width=\linewidth]{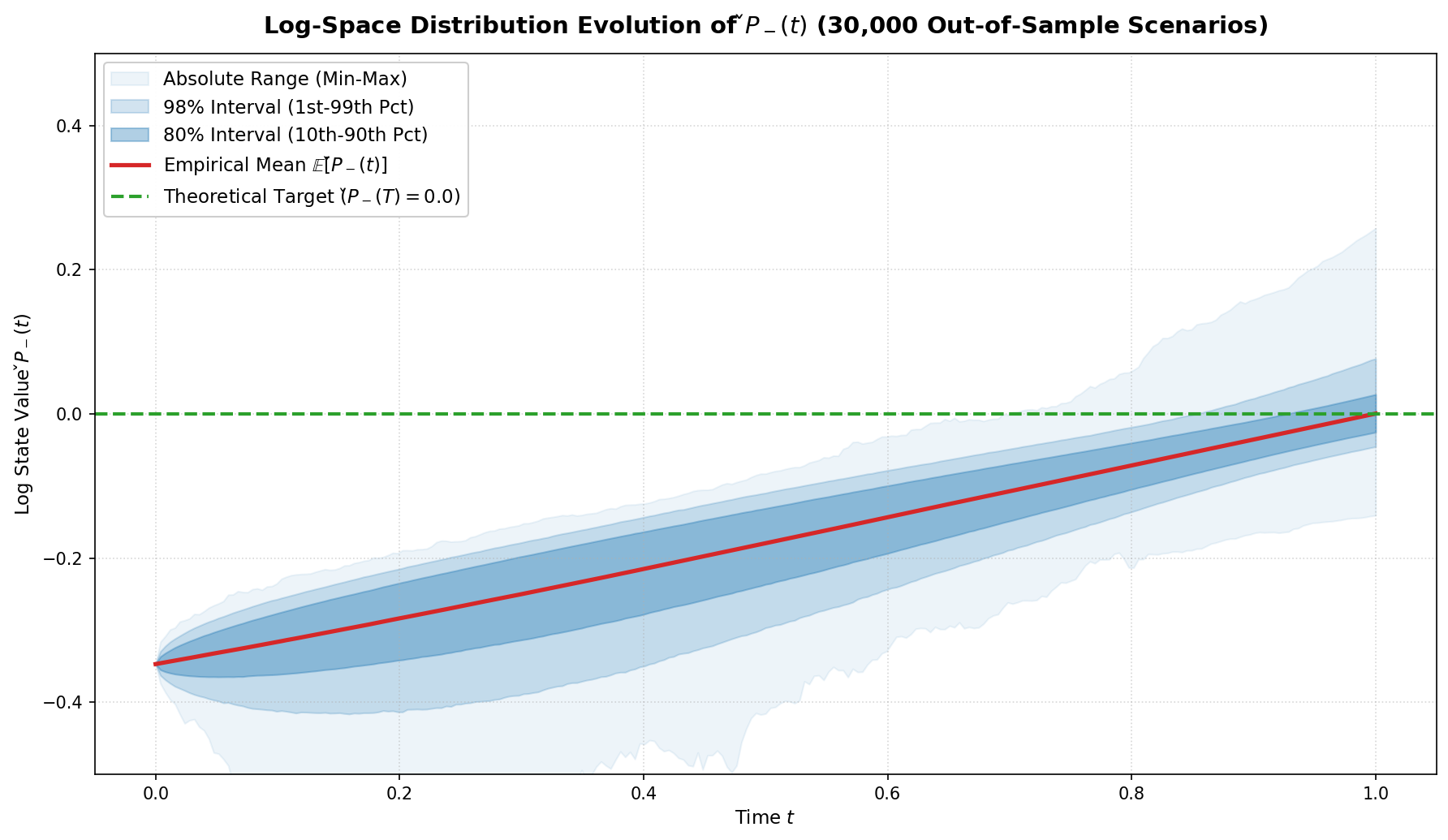}
        \caption*{$\check P_-$}
    \end{minipage}
    \caption{Out-of-sample evolution of $\check P_{\pm}(t)$ over
    $30{,}000$ independent Monte Carlo paths.}
    \label{fig:oos_funnel}
\end{figure}

We further evaluate the trained networks on $30{,}000$ independent Monte Carlo paths. Figure~\ref{fig:oos_funnel} shows the empirical mean of $\check P_{\pm}(t)$ together with the central $80\%$ and $98\%$ empirical
quantile bands and the full sample range. For both equations, the empirical means approach
zero, and the terminal distributions are
concentrated around the prescribed terminal
condition.

Table~\ref{tab:oos_metrics} summarizes the corresponding terminal
statistics. The terminal means and mean-square errors are small for both
equations, showing that the trained networks reproduce the prescribed
terminal condition well on independent simulated paths.

\begin{table}[htbp]
    \centering
    \caption{Out-of-sample terminal statistics over $30{,}000$ independent
    Monte Carlo paths. The theoretical terminal condition is
    $\check P_{\pm}(T)=0$.}
    \label{tab:oos_metrics}
    \renewcommand{\arraystretch}{1.15}
    \begin{tabular}{lcc}
        \hline
        \textbf{Metric}
        & \textbf{$\check P_+$}
        & \textbf{$\check P_-$} \\
        \hline
        Terminal mean
        & $0.000422$
        & $0.000469$ \\
        Terminal standard deviation
        & $0.023692$
        & $0.023564$ \\
        Terminal MSE
        & $5.6150\times10^{-4}$
        & $5.5548\times10^{-4}$ \\
        Maximum absolute deviation
        & $0.3678$
        & $0.6086$ \\
        \hline
    \end{tabular}
\end{table}

\section{Concluding Remarks}

This paper developed a reflected Markovian factor framework for constrained mean--variance portfolio selection with random coefficients. The reflected-factor dynamics, together with
continuous state-dependent market coefficients,
ensure boundedness of the market coefficients
and naturally induce homogeneous Neumann
boundary conditions. A logarithmic transformation reformulates the associated ESREs as quadratic BSDEs, whose Markovian representations are linked to semilinear Neumann problems. We established the existence and uniqueness of bounded continuous viscosity solutions to these PDEs. We further developed a deep BSDE method combining reflected-path simulation, geometry-aware boundary treatment, and FISTA-based evaluation of the constrained Hamiltonian, providing a mesh-free approach to high-dimensional ESREs.

Future research may extend the reflected-factor approach to other constrained portfolio optimization problems, including those with CARA and CRRA preferences. Another direction is to develop error estimates for the proposed deep BSDE method and investigate its application to more general high-dimensional PDEs with Neumann boundary conditions.

\appendix
\setcounter{theorem}{0}
\renewcommand{\thetheorem}{\thesection.\arabic{theorem}}
\setcounter{lemma}{0}
\renewcommand{\thelemma}{\thesection.\arabic{lemma}}
\setcounter{assumption}{0}
\renewcommand{\theassumption}{\thesection.\arabic{assumption}}

\section{Proof of the Kobylanski-type Transformation} \label{app-kobi}

\begin{proof}

We divide the proof into three steps.

\medskip
\noindent
\textbf{Step 1. Construction of the change of variables.}

Fix $J$ and $M$ as in the statement of the theorem. 
Under the transformation
\[
u(t,v)=e^{K(T-t)}\Psi(\widetilde u(t,v)),
\]
we have
\[
\Psi(\widetilde u(t,v))
=
e^{-K(T-t)}u(t,v)\in[-M,M],
\]
so it suffices to construct the change of variables on this bounded interval.

Let $C_H>1$ be a structural constant in
Proposition~\ref{prop:cone-generator}.
Set $C_0:=C_H+2C_H^2$ and choose $K>2C_0$.
Define $\Theta_K:=e^{KT}$. Choose $\lambda_K>8C_H\Theta_K$, and then choose $A_K$ sufficiently large so that
\[
A_K>
\max\left\{
2M\lambda_K+2,\,
\lambda_K+1
\right\}.
\]
For $y\in[-M,M]$, define
\[
w_K(y):=
A_K-\exp[-\lambda_K(y+M)].
\]
Then
\[
w_K(y)>0,
\qquad
w_K'(y)
=
\lambda_K\exp[-\lambda_K(y+M)]>0,
\qquad
w_K''(y)
=
-\lambda_K w_K'(y).
\]

Extend $w_K$ smoothly to $\mathbb R$ so that it remains bounded above
and below by strictly positive constants. Define $\Psi=\Psi_K$ as the
solution of
\[
\Psi'(\xi)=w_K(\Psi(\xi)),
\qquad
\Psi(0)=0.
\]
Since $w_K$ is smooth, bounded, and bounded away
from zero, the ODE admits a unique global solution.
Moreover, $\Psi$ is a strictly increasing 
$C^\infty$-diffeomorphism of $\mathbb R$ onto itself.
Then we set
\[
I_M:=\Psi^{-1}([-M,M]),
\qquad
p(\xi):=\Psi'(\xi),
\qquad
q(\xi):=\frac{\Psi''(\xi)}{\Psi'(\xi)}.
\]
For $\xi\in I_M$, writing $y=\Psi(\xi)$ gives
\[
p(\xi) = \Psi'(\xi) =w_K(y)>0,
\qquad
q(\xi)=w_K'(y)>0.
\]
Moreover,
\[
p'(\xi)=p(\xi)q(\xi), \qquad
q'(\xi)
=
w_K''(y)\Psi'(\xi)
=
-\lambda_K p(\xi)q(\xi),
\]
and
\[
\frac{q(\xi)}{p(\xi)}
=
\frac{w_K'(y)}{w_K(y)}
\leq
\frac{\lambda_K}{A_K-1}
\leq1.
\]

\medskip
\noindent
\textbf{Step 2. The transformed Neumann problem.}

Set $\theta(t):=e^{K(T-t)}$. For a smooth function $\widetilde u$, write $z:=\bm G(v)^*\nabla\widetilde u(t,v)$. Then
\[
\bm G(v)^*\nabla u
=
\theta(t)\Psi'(\widetilde u)z,
\]
and
\[
\partial_tu
=
\theta(t)\Psi'(\widetilde u)\partial_t\widetilde u
-
K\theta(t)\Psi(\widetilde u), \qquad
\mathcal Lu
=
\theta(t)\Psi'(\widetilde u)\mathcal L\widetilde u
+
\frac12\theta(t)\Psi''(\widetilde u)|z|^2.
\]

Substituting these identities into the original Neumann problem \eqref{eq:pde-sre} and dividing by
$\theta(t)\Psi'(\widetilde u)>0$ yields
\[
\partial_t\widetilde u
+
\mathcal L\widetilde u
+
\widetilde{\mathcal H}_K
\bigl(
t,v,\widetilde u,
\bm G(v)^*\nabla\widetilde u
\bigr)
=0,
\]
where
\begin{equation}\label{eq:transformed-H-app}
\begin{aligned}
\widetilde{\mathcal H}_K(t,v,\xi,z)
:={}&
-K\frac{\Psi(\xi)}{\Psi'(\xi)}
+
\frac12
\frac{\Psi''(\xi)}{\Psi'(\xi)}
|z|^2
+
\frac{1}{\theta(t)\Psi'(\xi)}
\mathcal H
\bigl(
t,v,\theta(t)\Psi'(\xi)z
\bigr).
\end{aligned}
\end{equation}

The terminal condition is preserved since $\theta(T)=1$ and
$\Psi(0)=0$:
\[
u(T,v)=0
\quad\Longleftrightarrow\quad
\widetilde u(T,v)=0.
\]
Likewise,
\[
\nabla u(t,v)
=
\theta(t)\Psi'(\widetilde u(t,v))
\nabla\widetilde u(t,v),
\]
and hence
\[
\langle\nabla u,\nabla\phi\rangle=0
\quad\Longleftrightarrow\quad
\langle\nabla\widetilde u,\nabla\phi\rangle=0.
\]
The corresponding terminal inequalities for viscosity
subsolutions and supersolutions are also preserved.

Since the transformation is smooth and strictly
increasing, applying the above calculations to
test functions shows that both the interior and
boundary operators are multiplied by the same
positive factor. Hence the transformation preserves
viscosity subsolutions and supersolutions.

\medskip
\noindent
\textbf{Step 3. Structural properties of the transformed Hamiltonian.}

Write
\[
\widetilde{\mathcal H}_K(t,v,\xi,z)
=
\frac12q(\xi)|z|^2
-
Kh(\xi)
+
\mathcal R(t,v,\xi,z),
\]
where
\[
h(\xi):=
\frac{\Psi(\xi)}{\Psi'(\xi)},
\qquad
\mathcal R(t,v,\xi,z)
:=
\frac{1}{\theta(t)p(\xi)}
\mathcal H
\bigl(
t,v,\theta(t)p(\xi)z
\bigr).
\]

We first establish strict monotonicity in $\xi$.
For fixed $(t,v,z)$, the map
\[
    \xi\mapsto\mathcal H(t,v,\theta(t)p(\xi)z)
\]
is Lipschitz continuous on $I_M$ by
Proposition~\ref{prop:cone-generator}, and hence
absolutely continuous and differentiable almost
everywhere. Moreover, the weighted Lipschitz estimate yields
\[
    \left|
    \frac{\mathrm d}{\mathrm d\xi}
    \mathcal H(t,v,\theta(t)p(\xi)z)
    \right|
    \leq
    C_H(1+2\theta(t)p(\xi)|z|)
    \theta(t)p'(\xi)|z|
\]
for a.e. $\xi\in I_M$.

Using the quotient rule, the identity $p'=pq$,
and the quadratic-growth estimate, we obtain
\begin{align*}
|\partial_\xi R(t,v,\xi,z)|
&\leq
\frac{q(\xi)}{\theta(t)p(\xi)}
|\mathcal H(t,v,\theta(t)p(\xi)z)| +
C_Hq(\xi)(1+2\theta(t)p(\xi)|z|)|z|\\
&\leq
C_H\frac{q(\xi)}{\theta(t)p(\xi)}
+C_Hq(\xi)|z|
+3C_H\theta(t)p(\xi)q(\xi)|z|^2
\end{align*}
for a.e. $\xi\in I_M$.

Since $q/p\leq1$ and $\theta(t)\geq1$,
Young's inequality gives
\[
    C_Hq(\xi)|z|
    \leq
    \frac{2C_H^2}{\lambda_K}
    +\frac{\lambda_K}{8}p(\xi)q(\xi)|z|^2.
\]
Consequently,
\[
    \partial_\xi R(t,v,\xi,z)
    \leq
    C_0+
    \left(3C_H\Theta_K+\frac{\lambda_K}{8}\right)
    p(\xi)q(\xi)|z|^2.
\]

% We first establish strict monotonicity in $\xi$.
% By Proposition~\ref{prop:cone-generator} and the smoothness of $p$,
% the map
% \[
% \xi\longmapsto\mathcal R(t,v,\xi,z)
% \]
% is locally Lipschitz on $I_M$, and is therefore differentiable for
% a.e.~$\xi\in I_M$ by Rademacher's theorem. At such points, the chain rule for Lipschitz functions gives
% \[
% \begin{aligned}
%     \partial_\xi\mathcal R(t,v,\xi,z) = -\frac{q(\xi)}{\theta(t)p(\xi)} \mathcal H\bigl(t,v,\theta(t)p(\xi)z\bigr) + q(\xi)
%     \left\langle
%         D_z\mathcal H\bigl(t,v,\theta(t)p(\xi)z\bigr), z
%     \right\rangle .
% \end{aligned}
% \]
% Since $q/p\leq1$, we have
% \[
%     C_H q(\xi)|z| \leq \frac{2C_H^2}{\lambda_K} + \frac{\lambda_K}{8}p(\xi)q(\xi)|z|^2.
% \]
% Notice that \(\lambda_K>8C_H\Theta_K\geq 1\), using 
% \[
%     |\mathcal H(t,v,z)|\leq C_H(1+|z|^2), \qquad
%     |D_z\mathcal H(t,v,z)|\leq C_H(1+|z|),
% \]
% at a.e. differentiability point, it follows that check there exists a constant $C_0>0$, independent of $K$, such that
% \[
% \partial_\xi\mathcal R(t,v,\xi,z)
% \leq
% C_0
% +
% \left(
% 2C_H\Theta_K+\frac{\lambda_K}{8}
% \right)
% p(\xi)q(\xi)|z|^2.
% \]
On the other hand,
\[
\frac{\mathrm d}{\mathrm d\xi}
\left(
\frac12q(\xi)|z|^2
\right)
=
-\frac{\lambda_K}{2}
p(\xi)q(\xi)|z|^2.
\]
Therefore,
\[
\partial_\xi
\left[
\frac12q(\xi)|z|^2
+
\mathcal R(t,v,\xi,z)
\right]
\leq
C_0
-
\left(
\frac{3\lambda_K}{8}
-
3C_H\Theta_K
\right)
p(\xi)q(\xi)|z|^2.
\]
Since $\lambda_K>8C_H\Theta_K$, the coefficient in parentheses is
strictly positive. Moreover, $p q$ is continuous and strictly positive
on the compact interval $I_M$. Thus, for some $\eta_K>0$,
\begin{equation}\label{eq:quadratic-absorption-app}
\partial_\xi
\left[
\frac12q(\xi)|z|^2
+
\mathcal R(t,v,\xi,z)
\right]
\leq
C_0-\eta_K|z|^2
\end{equation}
for a.e.~$\xi\in I_M$. It remains to control the zero-order term. A direct calculation yields
\[
h'(\xi)
=
1-\frac{y w_K'(y)}{w_K(y)}.
\]
Using $|y|\leq M$, $w_K(y)\geq A_K-1$, and
$w_K'(y)\leq\lambda_K$, the choice
$A_K>2M\lambda_K+2$ implies
\[
h'(\xi)\geq\frac12,
\qquad \xi\in I_M.
\]

Combining this estimate with
\eqref{eq:quadratic-absorption-app} gives
\[
\partial_\xi
\widetilde{\mathcal H}_K(t,v,\xi,z)
\leq
-\left(\frac K2-C_0\right)
-\eta_K|z|^2,
\]
for a.e.~$\xi\in I_M$. Setting
\[
    \kappa:=\min\left\{\frac K2-C_0,\eta_K\right\}>0,
\]
we obtain
\[
    \partial_\xi\widetilde{\mathcal H}_K(t,v,\xi,z)
    \leq-\kappa(1+|z|^2)
\]
for a.e. $\xi\in I_M$.

Therefore, for $\xi_1,\xi_2\in I_M$ with $\xi_1\geq\xi_2$,
\begin{equation*}
    \widetilde{\mathcal H}_K(t,v,\xi_1,z) -
\widetilde{\mathcal H}_K(t,v,\xi_2,z)=
\int_{\xi_2}^{\xi_1}
\partial_\xi
\widetilde{\mathcal H}_K(t,v,s,z)
\,\mathrm ds
\leq
-\kappa(\xi_1-\xi_2)(1+|z|^2).
\end{equation*}
This proves the desired strict finite-difference monotonicity.

Finally, once $K$ and $\Psi$ are fixed, the functions
\[
\Psi,\qquad
\Psi',\qquad
\frac1{\Psi'},\qquad
\frac{\Psi}{\Psi'},\qquad
\frac{\Psi''}{\Psi'},
\]
are bounded and Lipschitz continuous on $I_M$. Consequently,
Proposition~\ref{prop:cone-generator} implies that there exists
$C_{M,\Psi,K}>0$ such that
\[
|\widetilde{\mathcal H}_K(t,v,\xi,z)|
\leq
C_{M,\Psi,K}(1+|z|^2),
\]
and
\[
\begin{aligned}
& \quad |\widetilde{\mathcal H}_K(t_1,v_1,\xi_1,z_1)-
\widetilde{\mathcal H}_K(t_2,v_2,\xi_2,z_2)|
\\
&\leq
C_{M,\Psi,K}
(1+|z_1|^2+|z_2|^2)
\bigl(
|t_1-t_2|
+
|v_1-v_2|
+
|\xi_1-\xi_2|
\bigr)
+
C_{M,\Psi,K}
(1+|z_1|+|z_2|)
|z_1-z_2|.
\end{aligned}
\]
This completes the proof.
\end{proof}

\section{Proof of the Comparison Principle}
\label{app-comparison}

For a function $u:[0,T]\times\overline{\mathcal D}\to\mathbb R$,
we say that
\[
(a,p,X)\in\mathcal P^{2,+}u(t,x),
\]
if
\[
\begin{aligned}
u(s,y)\leq
u(t,x)+a(s-t)+\langle p,y-x\rangle+\frac12\langle X(y-x),y-x\rangle
+o\bigl(|s-t|+|y-x|^2\bigr),
\end{aligned}
\]
as $(s,y)\to(t,x)$ within $[0,T]\times\overline{\mathcal D}$. The parabolic subjet $\mathcal P^{2,-}u(t,x)$ is defined by reversing the inequality, and $\overline{\mathcal P}^{2,\pm}$ denote the corresponding sequential closures.

We use the shorthand
\[
\mathfrak N(p,v)
:=
-\langle p,\nabla\phi(v)\rangle,
\qquad
p\in\mathbb R^n,\quad v\in\partial\mathcal D.
\]
Thus, for a smooth test function $\psi$,
\[
\mathfrak N[\psi](t,v)
=
\mathfrak N(\nabla\psi(t,v),v).
\]

\begin{proof}[Proof of Theorem~\ref{thm:comparison-transformed-neumann}]

Suppose, by contradiction, that
\[
M_*
:=
\max_{[0,T]\times\overline{\mathcal D}}
\bigl(\widetilde u-\widetilde w\bigr)>0.
\]
Since
\[
\widetilde u(T,v)\leq \widetilde w(T,v),
\qquad v\in\overline{\mathcal D},
\]
a maximizer of $\widetilde u-\widetilde w$ corresponding to the positive value $M_*$ cannot occur at the terminal time.

Let $\alpha>0$ be a small constant, to be fixed below. For
$\varepsilon,\delta>0$, define
\begin{equation}\label{eq:double-penalty}
\begin{aligned}
\Phi_{\varepsilon,\delta,\alpha}(t,s,x,y)
:= \widetilde u(t,x)-\widetilde w(s,y)
-\frac{|x-y|^2}{2\varepsilon}
-\frac{|t-s|^2}{2\delta} +
\alpha\bigl(\phi(x)+\phi(y)\bigr),
\end{aligned}
\end{equation}
for
\[
(t,s,x,y)
\in
[0,T]^2\times
\overline{\mathcal D}^2.
\]
Since $\widetilde u$ is upper semicontinuous,
$\widetilde w$ is lower semicontinuous, and the underlying domain is
compact, $\Phi_{\varepsilon,\delta,\alpha}$ attains its maximum at some
point
\[
(t_{\varepsilon,\delta},
 s_{\varepsilon,\delta},
 x_{\varepsilon,\delta},
 y_{\varepsilon,\delta}).
\]
For simplicity, we write
\[
(t_\epsilon,s_\epsilon,x_\epsilon,y_\epsilon)
:=
(t_{\varepsilon,\delta},
 s_{\varepsilon,\delta},
 x_{\varepsilon,\delta},
 y_{\varepsilon,\delta}).
\]

\medskip
\noindent
\textbf{Step 1. Doubling-of-variables estimates.}

Define
\[
\widetilde u_\alpha(t,x)
:=
\widetilde u(t,x)+\alpha\phi(x),
\qquad
\widetilde w_\alpha(t,x)
:=
\widetilde w(t,x)-\alpha\phi(x).
\]
By the standard doubling-of-variables argument,

\[
|x_\epsilon-y_\epsilon|\to0,
\qquad
|t_\epsilon-s_\epsilon|\to0,
\qquad
\frac{|x_\epsilon-y_\epsilon|^2}{\varepsilon}
+
\frac{|t_\epsilon-s_\epsilon|^2}{\delta}
\longrightarrow0.
\]
as $\varepsilon,\delta\to0$, up to a subsequence. Moreover, for sufficiently small $\varepsilon,\delta$, we have
\begin{equation}\label{eq:positive-gap}
t_\epsilon<T,
\qquad
s_\epsilon<T,
\qquad
\widetilde u_\epsilon-\widetilde w_\epsilon
\geq\frac{M_*}{2},
\end{equation}
where
\[
\widetilde u_\epsilon
:=
\widetilde u(t_\epsilon,x_\epsilon),
\qquad
\widetilde w_\epsilon
:=
\widetilde w(s_\epsilon,y_\epsilon).
\]

\medskip
\noindent
\textbf{Step 2. Treatment of the Neumann boundary.}

We first recall a simple consequence of the geometry of the defining
function $\phi$. Since
\[
\phi(x)=0
\quad\text{for }x\in\partial\mathcal D,
\qquad
\phi(y)\geq0
\quad\text{for }y\in\overline{\mathcal D},
\]
Taylor's formula and the boundedness of $D^2\phi$ imply that there exists
$C_D>0$ such that
\begin{equation*}\label{eq:boundary-geometry}
\langle\nabla\phi(x),y-x\rangle
\geq
-C_D|x-y|^2,
\qquad
x\in\partial\mathcal D,\quad
y\in\overline{\mathcal D}.
\end{equation*}

Set
\[
a_\epsilon
:=
\frac{t_\epsilon-s_\epsilon}{\delta},
\qquad
p_\epsilon
:=
\frac{x_\epsilon-y_\epsilon}{\varepsilon},
\]
and
\begin{equation*}\label{eq:px-py}
p_\epsilon^x
:=
p_\epsilon-\alpha\nabla\phi(x_\epsilon),
\qquad
p_\epsilon^y
:=
p_\epsilon+\alpha\nabla\phi(y_\epsilon).
\end{equation*}

Suppose first that $x_\epsilon\in\partial\mathcal D$. Then
\[
\begin{aligned}
\mathfrak N(p_\epsilon^x,x_\epsilon)
&=
-\langle
p_\epsilon^x,\nabla\phi(x_\epsilon)
\rangle \geq
\alpha
-
C_D
\frac{|x_\epsilon-y_\epsilon|^2}{\varepsilon} > 0,
\end{aligned}
\]
for sufficiently small $\epsilon$, where we have used $|\nabla\phi|=1$ on $\partial\mathcal D$. Hence, at $x_\epsilon$, the relaxed viscosity subsolution condition $\min\{\mathfrak R,\mathfrak N\}\leq0$ cannot be satisfied through the boundary branch alone. Therefore the PDE branch must satisfy
\[
\mathfrak R\leq0.
\]
If $x_\epsilon\in\mathcal D$, the same PDE inequality follows directly from the interior viscosity condition.

Similarly, if $y_\epsilon\in\partial\mathcal D$, then
\[
\mathfrak N(p_\epsilon^y,y_\epsilon)
\leq -\alpha+
C_{\mathcal D}\frac{|x_\epsilon-y_\epsilon|^2}{\epsilon}
<0
\]
for sufficiently small $\epsilon,\delta$.
Hence the PDE branch of the viscosity supersolution
inequality must also hold.

\medskip
\noindent
\textbf{Step 3. Application of Crandall--Ishii theorem.}

Define the smooth penalization
\[
\begin{aligned}
\varphi_{\varepsilon,\delta,\alpha}(t,s,x,y)
:=
\frac{|t-s|^2}{2\delta}
+
\frac{|x-y|^2}{2\varepsilon}
-
\alpha\bigl(\phi(x)+\phi(y)\bigr).
\end{aligned}
\]
Then
\[
\Phi_{\varepsilon,\delta,\alpha}
=
\widetilde u(t,x)-\widetilde w(s,y)
-
\varphi_{\varepsilon,\delta,\alpha}(t,s,x,y).
\]

We apply the Crandall--Ishii theorem (Theorem~3.2 in \cite{crandall1992user}) of sums to the functions
$\widetilde u$ and $-\widetilde w$ on the compact set
$[0,T]\times\overline{\mathcal D}$, treating time
and space jointly as variables in $\mathbb R^{1+n}$.
We obtain symmetric matrices $\mathbb X_\epsilon,\mathbb Y_\epsilon \in\mathbb S^{n+1}$. Taking their spatial principal blocks $X_\epsilon,Y_\epsilon\in\mathbb S^n$
and passing from the elliptic to the parabolic
closed semijets, we obtain
% We regard $(t,x)$ and $(s,y)$ as variables in $\mathbb R^{1+n}$ and apply the Crandall--Ishii theorem (see, for example, Theorem 3.2 in Crandall et al. \cite{crandall1992user}) of sums to
% $\widetilde u(t,x)$ and $-\widetilde w(s,y)$ at the maximum point
% $(t_\epsilon,s_\epsilon,x_\epsilon,y_\epsilon)$.
% We obtain symmetric matrices $\mathbb X_\epsilon,\mathbb Y_\epsilon \in\mathbb S^{n+1}$ such that the corresponding jets have gradients $(a_\epsilon,p_\epsilon^x)$ and $(a_\epsilon,p_\epsilon^y)$ respectively. Denoting by $X_\epsilon,Y_\epsilon\in\mathbb S^n$ the spatial principal blocks of
% $\mathbb X_\epsilon,\mathbb Y_\epsilon$, we obtain
\begin{equation*}\label{eq:jets}
(a_\epsilon,p_\epsilon^x,X_\epsilon)
\in
\overline{\mathcal P}^{2,+}
\widetilde u(t_\epsilon,x_\epsilon),
\end{equation*}
and
\begin{equation*}\label{eq:jets2}
(a_\epsilon,p_\epsilon^y,Y_\epsilon)
\in
\overline{\mathcal P}^{2,-}
\widetilde w(s_\epsilon,y_\epsilon).
\end{equation*}
In particular, the temporal components of the two jets coincide.

The matrix inequality in the theorem of sums further implies, after
restricting to the spatial blocks, that
\begin{equation}\label{eq:spatial-ishii}
\begin{pmatrix}
X_\epsilon & 0\\
0 & -Y_\epsilon
\end{pmatrix}
\leq
\frac{C}{\varepsilon}
\begin{pmatrix}
I & -I\\
-I & I
\end{pmatrix}
+
C\alpha I,
\end{equation}
where $C$ is independent of $\varepsilon,\delta$ and $\alpha$. To derive the spatial matrix estimate \eqref{eq:spatial-ishii}, write
\[
K_I:=
\begin{pmatrix}
I&-I\\
-I&I
\end{pmatrix},
\qquad
D_\epsilon:=
\begin{pmatrix}
D^2\phi(x_\epsilon)&0\\
0&D^2\phi(y_\epsilon)
\end{pmatrix}.
\]
The spatial block of
$A_\epsilon=D^2\varphi_{\varepsilon,\delta,\alpha}$ is
\[
B_\epsilon
=
\frac1\varepsilon K_I-\alpha D_\epsilon.
\]
Since the time and space variables are uncoupled in the penalization,
the spatial principal block of the Ishii matrix inequality gives
\[
\begin{pmatrix}
X_\epsilon&0\\
0&-Y_\epsilon
\end{pmatrix}
\leq
B_\epsilon+\rho B_\epsilon^2.
\]
Choosing $\rho=\varepsilon$ and using $K_I^2=2K_I$, we obtain
\[
B_\epsilon+\varepsilon B_\epsilon^2
=
\frac{3}{\varepsilon}K_I
-\alpha D_\epsilon
-\alpha(K_ID_\epsilon+D_\epsilon K_I)
+\varepsilon\alpha^2D_\epsilon^2.
\]
Since $D^2\phi$ is bounded on $\overline{\mathcal D}$, the last three
terms are bounded above by $C\alpha I$ for
$0<\varepsilon,\alpha\leq1$. Hence
\[
\begin{pmatrix}
X_\epsilon&0\\
0&-Y_\epsilon
\end{pmatrix}
\leq
\frac{C}{\varepsilon}
\begin{pmatrix}
I&-I\\
-I&I
\end{pmatrix}
+C\alpha I.
\]

Consequently, using the Lipschitz continuity and boundedness of
$\bm G$,
\begin{equation}\label{eq:diffusion-estimate}
\begin{aligned}
\frac12
\operatorname{Tr}\Big(
\bm G(x_\epsilon)\bm G(x_\epsilon)^*X_\epsilon
-
\bm G(y_\epsilon)\bm G(y_\epsilon)^*Y_\epsilon
\Big)\leq
C
\frac{|x_\epsilon-y_\epsilon|^2}{\varepsilon}
+
C\alpha.
\end{aligned}
\end{equation}
% Indeed, multiplying \eqref{eq:spatial-ishii} by the block matrix
% constructed from $\bm G(x_\epsilon)$ and $\bm G(y_\epsilon)$ gives a
% term proportional to
% \[
% \frac{1}{\varepsilon}
% |\bm G(x_\epsilon)-\bm G(y_\epsilon)|^2,
% \]
% which is bounded by
% $C|x_\epsilon-y_\epsilon|^2/\varepsilon$.

Similarly, by the Lipschitz continuity and boundedness of $\bm F$,
\begin{equation}\label{eq:drift-estimate}
\begin{aligned}
\langle\bm F(x_\epsilon),p_\epsilon^x\rangle
-
\langle\bm F(y_\epsilon),p_\epsilon^y\rangle
\leq
C
\frac{|x_\epsilon-y_\epsilon|^2}{\varepsilon}
+
C\alpha.
\end{aligned}
\end{equation}
Combining \eqref{eq:diffusion-estimate} and
\eqref{eq:drift-estimate}, we obtain
\begin{equation}\label{eq:L-estimate}
\mathscr A_\epsilon
\leq
C
\frac{|x_\epsilon-y_\epsilon|^2}{\varepsilon}
+
C\alpha,
\end{equation}
where
\[
\begin{aligned}
\mathscr A_\epsilon
:=
\frac12
\operatorname{Tr}\Big(
\bm G(x_\epsilon)\bm G(x_\epsilon)^*X_\epsilon
-
\bm G(y_\epsilon)\bm G(y_\epsilon)^*Y_\epsilon
\Big)
+
\langle\bm F(x_\epsilon),p_\epsilon^x\rangle
-
\langle\bm F(y_\epsilon),p_\epsilon^y\rangle.
\end{aligned}
\]
\medskip
\noindent
\textbf{Step 4. Viscosity inequalities and contradiction.}

Set
\[
z_\epsilon^x
:=
\bm G(x_\epsilon)^*p_\epsilon^x,
\qquad
z_\epsilon^y
:=
\bm G(y_\epsilon)^*p_\epsilon^y.
\]
By continuity of the coefficients and the Hamiltonian,
the viscosity inequalities extend to the closed
parabolic semijets. By Step 2, the PDE branches
of these inequalities are active at both points.
Hence
\[
\begin{aligned}
&-a_\epsilon
-\frac12\operatorname{Tr}\bigl(G(x_\epsilon)\bm G(x_\epsilon)^*X_\epsilon\bigr)
-\langle\bm F(x_\epsilon),p_\epsilon^x\rangle
-\widetilde{\mathcal H}_K
\left(
t_\epsilon,x_\epsilon,
\widetilde u_\epsilon,
z_\epsilon^x
\right)
\leq0,
\end{aligned}
\]
whereas
\[
\begin{aligned}
&-a_\epsilon
-\frac12\operatorname{Tr}\bigl(G(y_\epsilon)\bm G(y_\epsilon)^* Y_\epsilon\bigr)
-\langle\bm F(y_\epsilon),p_\epsilon^y\rangle
-\widetilde{\mathcal H}_K
\left(
s_\epsilon,y_\epsilon,
\widetilde w_\epsilon,
z_\epsilon^y
\right)
\geq0.
\end{aligned}
\]
The two occurrences of $a_\epsilon$ cancel. Using the estimates obtained above for the
second- and first-order terms, we therefore obtain
\begin{equation}\label{eq:comparison-reduced}
0
\leq
C\frac{|x_\epsilon-y_\epsilon|^2}{\varepsilon}
+
C\alpha
+
\Delta\widetilde{\mathcal H}_\epsilon,
\end{equation}
where
\[
\begin{aligned}
\Delta\widetilde{\mathcal H}_\epsilon
:={}&
\widetilde{\mathcal H}_K
(t_\epsilon,x_\epsilon,
\widetilde u_\epsilon,z_\epsilon^x)
-
\widetilde{\mathcal H}_K
(s_\epsilon,y_\epsilon,
\widetilde w_\epsilon,z_\epsilon^y).
\end{aligned}
\]

We decompose
\[
\begin{aligned}
\Delta\widetilde{\mathcal H}_\epsilon
=
\widetilde{\mathcal H}_K
(t_\epsilon,x_\epsilon,
\widetilde u_\epsilon,z_\epsilon^x)
-
\widetilde{\mathcal H}_K
(t_\epsilon,x_\epsilon,
\widetilde w_\epsilon,z_\epsilon^x)
+
\widetilde{\mathcal H}_K
(t_\epsilon,x_\epsilon,
\widetilde w_\epsilon,z_\epsilon^x)
-
\widetilde{\mathcal H}_K
(s_\epsilon,y_\epsilon,
\widetilde w_\epsilon,z_\epsilon^y).
\end{aligned}
\]
By the strict monotonicity of
$\widetilde{\mathcal H}_K$ and
$\widetilde u_\epsilon-\widetilde w_\epsilon\geq M_*/2$,
\[
\widetilde{\mathcal H}_K
(t_\epsilon,x_\epsilon,
\widetilde u_\epsilon,z_\epsilon^x)
-
\widetilde{\mathcal H}_K
(t_\epsilon,x_\epsilon,
\widetilde w_\epsilon,z_\epsilon^x)
\leq
-\frac{\kappa M_*}{2}
\bigl(1+|z_\epsilon^x|^2\bigr).
\]
Moreover, it is easy to check
\[
|z_\epsilon^x-z_\epsilon^y|
\leq
C\frac{|x_\epsilon-y_\epsilon|^2}{\varepsilon}
+
C\alpha.
\]
Hence, by the weighted continuity estimate for
$\widetilde{\mathcal H}_K$ and Young's inequality, for every $\eta>0$,
\[
\begin{aligned}
& \quad
\widetilde{\mathcal H}_K
(t_\epsilon,x_\epsilon,
\widetilde w_\epsilon,z_\epsilon^x)
-
\widetilde{\mathcal H}_K
(s_\epsilon,y_\epsilon,
\widetilde w_\epsilon,z_\epsilon^y)
\\
&\leq
\Big[
C\omega_f(
|t_\epsilon-s_\epsilon|
+
|x_\epsilon-y_\epsilon|
)
+\eta
\Big]
\bigl(1+|z_\epsilon^x|^2\bigr)
+
C_\eta
\left(
\frac{|x_\epsilon-y_\epsilon|^4}{\varepsilon^2}
+\alpha^2
\right).
\end{aligned}
\]
Choosing $\eta=\kappa M_*/8$ and then
$\varepsilon,\delta$ sufficiently small, we obtain
\[
\Delta\widetilde{\mathcal H}_\epsilon
\leq
-\frac{\kappa M_*}{4}
\bigl(1+|z_\epsilon^x|^2\bigr)
+
C_{M_*}
\left(
\frac{|x_\epsilon-y_\epsilon|^4}{\varepsilon^2}
+\alpha^2
\right).
\]
Substituting this into \eqref{eq:comparison-reduced} gives
\[
0
\leq
-\frac{\kappa M_*}{4}
+
C\alpha
+
C_{M_*}\alpha^2
+
o_{\varepsilon,\delta}(1),
\]
where we used
\[
\frac{|x_\epsilon-y_\epsilon|^2}{\varepsilon}\to0.
\]
Finally, choose $\alpha>0$ sufficiently small and then let
$\varepsilon,\delta\to0$. We obtain
\[
0\leq-\frac{\kappa M_*}{8}<0,
\]
a contradiction. Therefore
\[
\widetilde u(t,v)\leq\widetilde w(t,v),
\qquad
(t,v)\in[0,T]\times\overline{\mathcal D}.
\]

\end{proof}

\bibliographystyle{amsplain}

\bibliography{reference}

\end{document}